\documentclass{amsart}

\usepackage{amssymb}
\usepackage{amsmath}
\usepackage{hyperref}
\usepackage{bm}
\usepackage[cal=boondoxo]{mathalfa}

\newcommand{\Gam}{\Gamma_{\mathbb{R}}}
\newcommand{\half}{\frac{1}{2}}
\newcommand{\thalf}{\tfrac{1}{2}}
\newcommand{\summ}{\mathop{{\sum}^{\star}}}

\newcommand{\sym}{{\rm sym}}
\newcommand{\intt}{\int_{-\infty}^{\infty}}

\newcommand{\jj}{\mathcal{j}}
\newcommand{\kk}{\mathcal{k}}

\numberwithin{equation}{section}

\newtheorem{theorem}{Theorem}[section]
\newtheorem{proposition}[theorem]{Proposition}
\newtheorem{lemma}[theorem]{Lemma}

\makeatletter
\@namedef{subjclassname@2010}{%
  \textup{2010} Mathematics Subject Classification}
\makeatother

\begin{document}

\title{On the fourth moment of Hecke Maass forms and the Random Wave Conjecture}

\author{Jack Buttcane}
\address{
Mathematics Department\\ 244 Mathematics Building\\ Buffalo, NY 14260 \\ USA}
\email{buttcane@buffalo.edu}
\thanks{The first author was supported by National Science Foundation grant DMS-1601919}

\author{Rizwanur Khan}
\address{
Science Program\\ Texas A\&M University at Qatar\\ PO Box 23874\\ Doha, Qatar}
\email{rizwanur.khan@qatar.tamu.edu }

\subjclass[2010]{Primary: 11F12, 11M99; Secondary: 81Q50}

\keywords{$L$-functions, automorphic forms, $L^4$-norm, quantum chaos, random wave conjecture}

\begin{abstract} Conditionally on the Generalized Lindel\"{o}f Hypothesis, we obtain an asymptotic for the fourth moment of Hecke Maass cusp forms of large Laplacian eigenvalue for the full modular group. This lends support to the Random Wave Conjecture.
 \end{abstract}

\maketitle

\section{Introduction}

A central question in Arithmetic Quantum Chaos concerns the distribution of Hecke Maass forms for the full modular group. The Random Wave Conjecture (RWC) predicts that for large Laplacian eigenvalue, the distribution is close to random. One way to formulate this is to conjecture that on fixed compact sets, the moments of Hecke Maass cusp forms of large Laplacian eigenvalue asymptotically equal the moments of a random variable with Gaussian distribution. Until now only numerical work and heuristic arguments (see \cite{hejrac, hejstr}) have supported this conjecture. On the theoretical side, the fourth moment in particular is a natural and important case to study, as it is reduces, via Watson's formula, to a problem on $L$-functions. But proving an asymptotic for the fourth moment of Hecke Maass forms seems to be beyond the reach of current technology. The goal of this paper is to establish such an asymptotic (over $\Gamma\backslash \mathbb{H}$) on the assumption of the Generalized Lindel\"{o}f Hypothesis (GLH). One may wonder what the benefit is of proving one conjecture based on another. One answer is that these conjectures are unrelated. Another answer is of course that the GLH is a much more well accepted conjecture in mathematics and its truth is very firmly believed. This puts the RWC on more solid ground. Our main result is the following.

\begin{theorem} \label{main} Assume the GLH. Let $f$ be an even or odd Hecke-Maass cusp form for $ \Gamma= \text{SL}_2(\mathbb{Z})$ with Laplacian eigenvalue $\lambda=\frac{1}{4}+ T^2$, where $T>0$. Let $f$ be normalized as follows:
\begin{align}
\label{probmeas} \frac{1}{\int_{ \Gamma \backslash \mathbb{H}} 1 \frac{dxdy}{y^2} } \int_{ \Gamma \backslash \mathbb{H}} f(z)^2 \  \frac{dxdy}{y^2}=1.
\end{align}
There exists a constant $\delta>0$ such that
\begin{align}
\label{rwc} \frac{1}{\int_{ \Gamma \backslash \mathbb{H}} 1 \frac{dxdy}{y^2}} \int_{ \Gamma \backslash \mathbb{H}} f(z)^4 \  \frac{dxdy}{y^2} =  \frac{1}{\sqrt{2\pi}}\intt t^{4} e^{\frac{-t^2}{2}} dt +O(T^{-\delta})
\end{align}
as $T \to \infty$. 
\end{theorem}
\noindent Thus our result confirms on GLH a prediction of the RWC, with a power saving. Note that $f$ is real valued because it is assumed to be even or odd (and of weight 0, with trivial nebentypus).  The left hand side of (\ref{rwc}) will be studied by first using Watson's formula to relate it to a mean value of $L$-functions. If the GLH were applied at this point, it would immediately yield the upper bound $O(T^\epsilon)$. To go beyond this and obtain an asymptotic with power saving, even on GLH, requires a lot of work and the full power of spectral theory.  Some care is also needed to avoid reliance on the Ramanujan Conjecture (both at the finite and infinite places). Our proof will show that if not GLH, then at the very least what is required is a subconvexity bound for the $L$-function associated to the Rankin-Selberg product of an (essentially) fixed Hecke cusp form and the symmetric square lift of $f$, in the eigenvalue aspect of $f$. See the discussion following Lemma \ref{wlem} and the last sentence of this paper. This GL(2)$\times$GL(3) subconvexity problem seems to be very difficult and, interestingly, is also essentially what is required to get a power saving error in the QUE problem. It is safe to say that our problem is more difficult than the quantitative QUE problem.

There have been some unconditional results short of an asymptotic for the fourth moment of automorphic forms on $\Gamma\backslash \mathbb{H}$. In the eigenvalue aspect, an essentially optimal upper bound was proven by Spinu \cite{spi} for the fourth moment of truncated Eisenstein series, and by Luo \cite{luo} for dihedral Maass forms. Both of these results hinged on the spectral large sieve. For Hecke Maass forms, Sarnak and Watson announced a sharp upper bound for the fourth moment in \cite[Theorem 3]{sar}, but a proof of this has not appeared. Holomorphic Hecke cusp forms of large weight are expected to be modelled by a complex Gaussian distribution (see \cite[Conjecture 1.2]{blokhayou}). However proving an asymptotic for the fourth moment in this case, on GLH or any other reasonable hypothesis, seems to be much harder than the problem for Maass forms. This is because the corresponding mean value of $L$-functions has a larger ``log of conductor to log of family size ratio'' (see the discussion below). The best known upper bound for the fourth moment in the weight aspect is far from optimal; see \cite{blokhayou}.

Let $\{ u_j : j\ge 1 \}$ denote an orthonormal basis of even and odd Hecke Maass cusp forms for $\Gamma$, ordered by Laplacian eigenvalue $\frac{1}{4}+t_j^2$, where $t_j>1$.  Let $L(s,u_j)$ be the $L$-function attached to $u_j$, normalized so that its functional equation relates values at $s$ and $1-s$. The shape of the mean value of $L$-functions that we will need to evaluate is essentially 
\begin{align}
\label{shape} \sum_{t_j<2T} \frac{1}{t_j  T^\half (1+2T-t_j)^\half} \frac{L(\half,u_j) L(\half,u_j \times \sym^2 f)}{L(1,\sym^2 u_j)}.
\end{align}
The analytic conductors of $L(\half,u_j)$ and $ L(\half,u_j \times \sym^2 f)$ are $t_j^2$ and $t_j^2(1+|4T^2-t_j^2|)^2$ respectively. Thus the denominator above is about the same size as the convexity bound for the numerator. In the ``bulk'' range $T^{1-\epsilon}<t_j< 2T-T^{1-\epsilon}$, which is nearly a dyadic interval, the analytic conductor of the triple product $L$-function  $L(\half,u_j) L(\half,u_j \times \sym^2 f)$ is about $T^8$, while the sum is over about $T^2$ forms. Thus the ``log of conductor to log of family size ratio'' is 4. Our main work will be on treating this bulk range. The remaining ranges will handled immediately on the GLH. 

A similar mean value of triple product $L$-functions was considered (unconditionally) by Li in \cite{li-x}, but there the $GL(3)$ form was fixed, while here it is not ($T$ tends to infinity). Thus our problem is clearly more complex. A similar mean value was also considered by the authors in \cite{butkha}, with $L(\half,u_j \times \sym^2 f)$ replaced by $L(\half,u_j \times \chi) L(u_j \times f)$, where $\chi$ is a quadratic character. Such a factorization occurs when $f$ is a dihedral form and this was the motivation for the work in \cite{butkha}, although in that paper we were not able to make any direct conclusions about the fourth moment. In this paper we use the methods of \cite{butkha} together with GL(3) Voronoi summation as a new ingredient to treat the present case where $L(\half,u_j \times \sym^2 f)$ does not factorize. The present analysis is more delicate, with the ranges of various parameters harder to control (in the same way that many problems in number theory involving the divisor function become more difficult when the divisor function is replaced by Fourier coefficients of cusp forms). For this reason that it is not clear a priori that our previous methods would work for this problem. Also, \cite{butkha} used a simplified weight function in place of the one given in (\ref{htdef}).

\section{Reduction to $L$-functions} \label{redsec}

\noindent {\bf Convention}. Throughout, $\epsilon$ will denote an arbitrarily small positive constant, which may not be the same from one occurrence to another. 

To prove Theorem \ref{main}, the starting point is to express the fourth moment as a mean value of $L$-functions. Note that since $\int_{ \Gamma \backslash \mathbb{H}} 1 \frac{dxdy}{y^2}=\frac{\pi}{3}$ and $\frac{1}{\sqrt{2\pi}}\intt t^{4} e^{\frac{-t^2}{2}} dt =3$, we need to show that with the normalization $\langle f, f\rangle = \langle f^2, 1\rangle = \frac{\pi}{3}$, where the inner product is the Petersson inner product, we have
\begin{align*}
\langle f^2, f^2 \rangle = \pi + O(T^{-\delta}).
\end{align*}
By the spectral theorem (see \cite[Theorem 15.5]{iwakow}) and Parseval's theorem, we have
\begin{align}
\nonumber \langle f^2, f^2 \rangle =  | \langle f^2,(\tfrac{3}{\pi})^{\half} \rangle|^2 +\frac{1}{4\pi} \int_{-\infty}^{\infty}  | \langle  f^2 , E(\cdot,  \thalf + it) \rangle|^2 \ dt +  \sum_{j\ge 1} | \langle  f^2 , u_j \rangle|^2 ,
\end{align}
where $E(z,s)$ is the standard Eisenstein series. By normalization of $f$, we have
\begin{align*}
| \langle f^2,(\tfrac{3}{\pi})^{\half} \rangle|^2= \frac{\pi}{3},
\end{align*}
and we will show that
\begin{lemma} \label{eislem} Let $f$ be as in Theorem \ref{main}. On the GLH, we have 
\begin{align}
\label{contpart} \int_{-\infty}^{\infty}  | \langle  f^2 , E(\cdot,  \thalf + it) \rangle|^2 \ dt \ll T^{-1+\epsilon}.
\end{align}
\end{lemma}
\noindent Thus it remains to prove
\begin{proposition} \label{prop} Let $f$ be as in Theorem \ref{main}. On the GLH, there exists some $\delta>0$ such that
\begin{align}
\label{remains} \sum_{j\ge 1} | \langle  f^2 , u_j \rangle|^2  = \frac{2\pi }{3}+O(T^{-\delta}).
\end{align}
\end{proposition}
\noindent Since $\sqrt{\tfrac{3}{\pi}} f$ has $L^2$-norm equal to 1, applying Watson's formula \cite[Theorem 3]{wat} (see also \cite[page 2624]{blokhayou}) to the inner product of $\frac{3}{\pi} f^2$ and $u_j$, we get
\begin{align*}
 | \langle \frac{3}{\pi} f^2 , u_j \rangle |^2 =  \frac{\pi}{2^3}H(t_j)\frac{L(\half,f\times f\times \overline{u_j})}{L(1,\sym^2 f)^2L(1,\sym^2 u_j)},
\end{align*} 
where the $L$-functions appearing above are defined in the next section and 
\begin{align}
\label{htdef} H(t) = \frac{ |\Gamma(\frac{\half+2iT+it}{2})|^2 |\Gamma(\frac{\half+2iT-it}{2})|^2 |\Gamma(\frac{\half+it}{2})|^4}{|\Gamma(\frac{1+2iT}{2})|^4  |\Gamma(\frac{1+2it}{2})|^2 }.
\end{align}
Equivalently, as $u_j$ is real, 
\begin{align}
\label{wats} | \langle f^2 , u_j \rangle |^2 =  \frac{\pi^3}{72}H(t_j)\frac{L(\half,u_j) L(\half, u_j \times \sym^2 f )}{L(1,\sym^2 f)^2L(1,\sym^2 u_j)}.
\end{align} 
For the weight function $H(t_j)$, we have by Stirling's approximation (see (\ref{stirs}) or \cite[Section 5.1.1]{spi}) that
\begin{align}
\label{stirht} H(t) = \frac{8\pi \exp\big(-\pi q(t,T)\big)}{ (1+|t|)\prod_\pm (1+|2T\pm t|)^{\half}}\Big\{1+O\Big(\frac{1}{1+|t|}+\frac{1}{1+|2T+t|}+\frac{1}{1+|2T-t|} \Big)\Big\},
\end{align}
where
\begin{align}
\label{qtT} q(t,T)=\begin{cases}
0 &\text{ for } |t|\le 2T\\
|t|-2T &\text{ for } |t|>2T.
\end{cases}
\end{align}
Thus the right hand side of (\ref{wats}) looks essentially like (\ref{shape}). 

\section{Preliminaries}

\subsection{Stirling's approximation}

For $\sigma>0$ fixed, as a first order approximation we have 
\begin{align}
\nonumber  &\Gamma(\sigma+i\gamma) = \sqrt{2\pi}  |\sigma+i\gamma|^{\sigma-\half+i\gamma} \exp\big( -\tfrac{\pi}{2} |\gamma| + i \text{sgn}(\gamma) ( \tfrac{\pi}{2}(\sigma-\thalf)- \gamma)\big)\Big(1 +O((1+|\gamma|)^{-1}) \Big),\\
\label{stirs} &|\Gamma(\sigma+i\gamma)| = \sqrt{2\pi}  (1+|\gamma|)^{\sigma-\half} \exp\big( -\tfrac{\pi}{2} |\gamma|\big)\Big(1 +O((1+|\gamma|)^{-1}) \Big)
\end{align}
where sgn$(\gamma)$ is $1$ if $\gamma$ is positive and $-1$ if $\gamma$ is negative. As $|\gamma|\to\infty$, this gives
\begin{align*}
\Gamma(\sigma+i\gamma) = \sqrt{2\pi}  |\gamma|^{\sigma-\half+i\gamma} \exp\big( -\tfrac{\pi}{2}|\gamma| + i \text{sgn}(\gamma) ( \tfrac{\pi}{2}(\sigma-\thalf)- \gamma)\big) \Big(1 +O(|\gamma|^{-1}) \Big).
\end{align*}

\subsection{Approximate functional equations}

Let $\lambda_j(n)$ and $\lambda_f(n)$ denote the (real) eigenvalues of the $n$-th Hecke operator corresponding to $u_j$ and $f$ respectively, where we write $\lambda_j(-n)=\lambda_j(n)$ for $u_j$ even and $\lambda_j(-n)=-\lambda_j(n)$ for $u_j$ odd. 
The $L$-function associated to $u_j$ is given by 
\begin{align*}
L(s, u_j) = \sum_{n\ge 1} \frac{\lambda_j(n)}{n^s}
  \end{align*}
 for $\Re(s)>1$.  
Let $A_f(n,1)=A_f(1,n)$ be given by
\begin{align*}
\sum_{n\ge 1} \frac{A_f(n,1)}{n^s} := \zeta(2s) \sum_{n\ge 1} \frac{\lambda_f(n^2)}{n^s}.
  \end{align*}
The right hand above equals $L(s, \sym^2 f)$ for $\Re(s)>1$. Now define $A_f(n,m)=A_f(m,n)$ by the Hecke relations
\begin{align}
\label{gl3hecke} A_f(n,m) = \sum_{v|(n,m)} \mu(v) A_f\Big(\frac{n}{v},1\Big) A_f\Big(\frac{m}{v},1\Big).
\end{align}  
With this, we can define
\begin{align*}
L(s, u_j \times \sym^2 f) = \sum_{m,r\ge 1} \frac{\lambda_j(m)A_f(m,r)}{(r^2m)^s}
\end{align*}
for $\Re(s)>1$.

Kim and Sarnak \cite[Appendix 2]{kimsar} have proven the following bounds towards the Ramanujan Conjecture:
\begin{align*}
&|\lambda_j(n)| \ll n^{\frac{7}{64}+\epsilon}\\
&|A_f(n,1)| \ll n^{\frac{7}{32}+\epsilon}.
\end{align*}
On average, the best possible bounds are known by \cite[Lemma 1]{iwa} and \cite[Lemma 2.2]{hofloc}:
\begin{align}
\nonumber &\sum_{n\le x} |\lambda_j(n)|^2 \ll x(|t_j|x)^\epsilon\\
\label{avgram2} &\sum_{n\le x} |A(n,1)|^2 \ll x(Tx)^\epsilon.
\end{align}
This implies by the Cauchy-Schwarz inequality and the Hecke relations (\ref{gl3hecke}) that
\begin{align}
\label{avgram} \sum_{n\le x} |\lambda_j(n)| \ll x(|t_j|x)^\epsilon,\\
\nonumber \sum_{\substack{n\le x\\m\le y}} |A(n,m)| \ll xy(Txy)^\epsilon.
\end{align}

Let 
\begin{align*}
&\Gam(s)= \pi^{-\frac{s}{2}} \Gamma(\tfrac{s}{2}),\\ 
&G_1(s)=\prod_\pm \Gam(s \pm it_j ),\\
&G_2(s)=\prod_\pm  \Gam(s \pm it_j +i2T )\Gam(s \pm it_j )  \Gam(s \pm it_j- i2T).
\end{align*}
For $u_j$ even we have the functional equations
\begin{align*}
&L(s, u_j)G_1(s)=L(1-s, u_j) G_1(1-s),\\
& L(s, u_j\times \sym^2 f)G_2(s)=L(1-s, u_j\times \sym^2 f)G_2(1-s).
\end{align*}
For $u_j$ odd we have the functional equation
\begin{align*}
 L(s, u_j)G_1(1+s)= - L(1-s, u_j) G_1(2-s),
 \end{align*}
 which implies that $L(\half, u_j)=0$.
All of these may be found in \cite[chapter 3]{gol} and \cite[page 1670]{li-x}.

We now set up approximate functional equations for the central values. But first we explain what we will need. Usually one takes an approximate functional equation with the shortest possible Dirichlet series. For $u_j$ even, this means taking
\begin{align*}
L(\thalf,u_j)\sim \sum_{n<t_j^{1+\epsilon}} \frac{\lambda_j(n)}{n^\half} + \sum_{n<t_j^{1+\epsilon}} \frac{\lambda_j(n)}{n^\half},
\end{align*}
and for $u_j$ odd, we could take
\begin{align*}
L(\thalf,u_j)\sim \sum_{n<t_j^{1+\epsilon}} \frac{\lambda_j(n)}{n^\half} - \sum_{n<t_j^{1+\epsilon}} \frac{\lambda_j(n)}{n^\half},
\end{align*}
which vanishes. Since $\lambda_j(-n)=\lambda_j(n)$ for $u_j$ even and $\lambda_j(-n)=-\lambda_j(n)$ for $u_j$ odd, in both cases we have
\begin{align*}
L(\thalf,u_j)\sim \sum_{n<t_j^{1+\epsilon}} \frac{\lambda_j(n)}{n^\half} + \sum_{n<t_j^{1+\epsilon}} \frac{\lambda_j(-n)}{n^\half}.
\end{align*}
As we will see below, to understand a mean value of the form $\sum_{T<t_j<2T} \lambda_j(\pm n)\lambda_j(m) $ using Kuznetsov's formula, the same sign Kuznetsov (the $+$ sign) leads to a $J$-Bessel transform while the opposite sign Kuznetsov (the $-$ sign) leads to a $K$-Bessel transform. Both transforms can be evaluated asymptotically and it turns out the main term of the $K$-Bessel transform has no oscillation. We find this easier to work with, so we reduce the analysis involving the same sign terms by taking an approximate functional equation with two Dirichlet series of unequal length as follows:
\begin{align*}
L(\thalf,u_j) \sim \sum_{n<t_j^{1+\epsilon}T^{-\beta}} \frac{\lambda_j(n)}{n^\half} + \sum_{n<t_j^{1+\epsilon}T^{\beta}} \frac{\lambda_j(-n)}{n^\half}.
\end{align*}
In this way, the Dirichlet series leading to the same sign terms is shorter. We now state this precisely.
\begin{lemma}\label{afelem}
For $u_j$ even, we have
\begin{align}
\label{afe1} &L(\thalf, u_j)=2\sum_{n\ge 1} \frac{\lambda_j(n)}{n^\half}V_1(n,t_j),\\
\label{afe2} &L(\thalf, u_j\times \sym^2 f)=2\sum_{m,r\ge 1} \frac{\lambda_j(m)A_f(r,m)}{rm^\half}V_2(r^2m,t_j),
\end{align}
where 
\begin{align}
\label{videf} V_i(x,t)=\frac{1}{2\pi i}\int_{(\sigma)} e^{s^2} x^{-s} \frac{G_i(\half+s)}{G_i(\half)} \frac{ds}{s}
\end{align}
for any $\sigma>0$.

Let $0<\alpha, \beta<\frac{1}{100}$ be some fixed constants to be determined later. For $u_j$ even or odd, and $T^{1-\alpha}<|t_j|<T^{1+\epsilon}$, we have
\begin{align}
\label{trickafe} L(\thalf, u_j)=\sum_\pm \sum_{n\ge 1} \frac{\lambda_j(\pm n)}{n^\half}V_1^\pm(n,t_j)+O(T^{-\half+\frac{\beta}{2}+\alpha+\epsilon}),
\end{align}
where
\begin{align*}
V_1^\pm(x,t)=\frac{1}{2\pi i}\int_{(\sigma)} e^{s^2}  (xT^{\pm \beta})^{-s} \frac{G_1(\half+s)}{G_1(\half)} \frac{ds}{s}
\end{align*}
for any $\sigma>0$.
\end{lemma}
\proof
This follows from \cite[Theorem 5.3]{iwakow} and the functional equations given above. For (\ref{afe1}) and (\ref{afe2}), take $G(u)=e^{u^2}$ and $X=1$ in that theorem. For (\ref{trickafe}), take $G(u)=e^{u^2}$ and $X=T^\beta$ to get
\begin{align}
\label{uba} L(\thalf, u_j)  = \sum_\pm\sum_{n\ge 1} \frac{\lambda_j( \pm n)}{n^{\half}}  \frac{1}{2\pi i} \int_{(\sigma)} e^{s^2} \left( x T^{\pm \beta} \right)^{-s}  \frac{G_1(\half+\kappa_j + s)}{G_1(\half +\kappa_k)}  \frac{ds}{s},
\end{align}
where $\kappa_j=0$ or $1$ as $u_j$ is even or odd. By the rapid decay of $e^{{s^2}}$ in vertical lines, we may restrict the integral above to $|\Im(s)|<T^{\epsilon}$. By Stirling's approximation, for $\Re(s)>0$ fixed, $|\Im(s)|<T^{\epsilon}$ and $|t|>T^{1-\alpha}$, we have
\begin{align}
\frac{G_1(\half+1+ s)}{G_1(\half +1)}= \frac{G_1(\half + s)}{G_1(\half )} + O(T^{-1+\alpha+\epsilon}).
\end{align}
 Thus up to a small error, the ratio of Gamma functions in (\ref{uba}) does not depend on $\kappa_j$. Also note that the sum in (\ref{uba}) can be restricted to $n<T^{1+\beta+\epsilon}$ up to admissible error by Stirling's approximation. Thus (\ref{trickafe}) follows, using (\ref{avgram}).
\endproof

Consider those values of $|t|$ that are roughly of size $2T$ but not too close to $2T$. That is, suppose that for some $0<\alpha<\frac{1}{100}$ to be fixed later, we have
\begin{align}
\label{closerange} T^{1-\alpha}< |t|<2T-T^{1-\alpha}.
\end{align}
In the integrals appearing in Lemma \ref{afelem}, write $s=\sigma+i\gamma$. By the rapid decay of $e^{s^2}$ in vertical lines, we may restrict these integrals to $|\gamma|<T^\epsilon$. We have by Stirling's approximation that
\begin{align}
\label{stirv1} V_1^\pm(x,t) = \frac{1}{2\pi i} \int_{(\sigma)} e^{s^2} \left( \frac{2\pi x T^{\pm \beta}   }{|t|} \right)^{-s} \Big( 1 + \sum_{n\le N} \frac{B_n(\sigma,\gamma)}{|t|^{n}} \Big)\frac{ds}{s} +O(T^{-\frac{N}{2}})
\end{align}
and
\begin{multline}
\label{stirv2} V_2(x,t) = \frac{1}{2\pi i} \int_{(\sigma)} e^{s^2} \left( \frac{8\pi^3  x }{|t(4T^2-t^2)|} \right)^{-s}  \Big( 1 + \sum_{n\le N}  \frac{B_n(\sigma,\gamma)}{|t|^{n}}\Big) \Big(1+ \sum_{n\le N} \frac{C_n(\sigma,\gamma)}{|2T+t|^{n}}\Big)\\ \times  \Big(1+ \sum_{n\le N} \frac{C_n(\sigma,\gamma)}{|2T-t|^{n}}\Big) \frac{ds}{s}+ O(T^{-\frac{N}{2}}),
\end{multline}
for any $N\ge 1$ and some $B_n(\sigma,\gamma)$ and $C_n(\sigma,\gamma)$ polynomial in $\sigma$ and $\gamma$. By thinking of $|t|$ and $|t(4T^2-t^2)|$ as being of size about $T$ and $T^3$ in the range (\ref{closerange}), and taking $\sigma$ as large as we like in the expressions (\ref{stirv1}-\ref{stirv2}), we see that the sums in (\ref{trickafe}) and (\ref{afe2}) have length about $T^{1\mp \beta}$ and $T^3$ respectively. We have
\begin{align}
\label{vder1} \frac{\partial^k}{\partial t^k} V_1^\pm(x,t)\ll |t|^{-k+\epsilon}\ll  T^{k(-1+\alpha+\epsilon)},
\end{align}
in the range (\ref{closerange}) by taking $N$ large enough and $\sigma=\epsilon$. Similarly,
\begin{align}
\label{vder2} \frac{\partial^k}{\partial t^k} V_2(x,t) \ll |t(4T^2-t^2)|^{-k+\epsilon}\ll T^{k(-3+2\alpha+\epsilon)} 
\end{align}
for $k \ge 0$.

\subsection{Kuznetsov trace formula}

Let
\begin{align*}
\lambda(n,t)=\sum_{ab=n} \Big(\frac{a}{b}\Big)^{it}
\end{align*}
and
\begin{align*}
\mathcal{J}^+(x, t) = \frac{2i}{\sinh(\pi t)} J_{2it}(4\pi x), \quad \mathcal{J}^-(x, t) = \frac{4}{\pi} K_{2it}(4\pi x) \cosh(\pi t).
\end{align*}
We have
\begin{lemma} \label{lem:Kuznetsov} 
Let $h(z)$ be an even, holomorphic function on $|\Im(z)| < \frac{1}{4}+\theta$ with decay $|h(z)| \ll (1+|z|)^{-2-\theta}$ on that strip, for some $\theta>0$. Then for $n,m > 0$, we have
\begin{multline}
\sum_{j \ge 1} \frac{\lambda_j(\pm n) \lambda_j(m)}{L(1,{\normalfont \text{sym}} ^2 u_j)} h(t_j) +\int_{-\infty}^\infty \frac{\lambda(n, t) \lambda(m, -t)}{|\zeta(1+2it)|^2} h(t) \frac{dt}{2\pi} \\
= \delta_{\pm n,m} \int_{-\infty}^\infty h(t) \frac{d^*t}{2\pi^2} + \sum_{c \ge 1} \frac{S(\pm n,m, c)}{c} \int_{-\infty}^\infty \mathcal{J}^\pm(\tfrac{\sqrt{nm}}{c}, t) h(t) \frac{d^* t}{2 \pi},
\end{multline}
where $\delta_{n,m}$ is 1 if $n=m$ and 0 otherwise (thus $\delta_{-n,m}$ is always 0) and $d^* t =  \tanh(\pi t) \, t dt$.
\end{lemma}
\proof

See \cite[Theorems 2.2, 2.4]{moto}. There, the function $h(z)$ must be holomorphic function on $|\Im(z)| < \frac{1}{2}+\theta$. The relaxation of this condition to $|\Im(z)| < \frac{1}{4}+\theta$ is due to Yoshida \cite{yos}. We need this version because $H(t)$, which was defined in section \ref{redsec}, has a pole at $t=\half i$.
\endproof 

For the Kuznetsov trace formula, we will need to asymptotically evaluate some Bessel transforms. The following lemmas, taken from \cite[Section 3.8]{butkha} together with the correction noted after Lemma \ref{kbeslemma} below, are analogous to the averages of real Bessel functions given in \cite[Corollary 8.2]{ils}:
\begin{align*}
&\sum_{k\equiv 0 \bmod 2} i^k J_{k-1}(x) h\Big(\frac{k}{K}\Big) \sim \text{oscillatory function supported on } x\gg K^{2-\epsilon}, \\
&\sum_{k\equiv 0 \bmod 2} J_{k-1}(x) h\Big(\frac{k}{K}\Big) \sim \text{non-oscillatory function supported on } x\asymp K.
\end{align*}

\begin{lemma} \label{jbeslemma} Let $0<\alpha<\frac{1}{100}$. For any $x>0$ and any smooth even  function $h$, compactly supported on $(T^{-\alpha}, T^{\alpha})\cup (-T^{\alpha}, -T^{-\alpha})$ with derivatives satisfying $\| h^{(k)} \|_\infty \ll (T^{\alpha})^{k}$, we have
\begin{multline}
\label{jbeslem} \intt \frac{J_{2it}(2\pi x)}{\cosh(\pi t)} h\Big(\frac{t}{T}\Big) t dt \\= \frac{-i\sqrt{2}}{\pi} \frac{T^2}{\sqrt{x}} \Re\left((1+i)e(x)   \int_0^\infty t h(t) e\Big(\frac{-t^2 T^2}{2\pi^2 x}\Big) dt \right) + O\Big( \frac{x}{T^{3-12\alpha}} \Big)+O(T^{-100}).
\end{multline}
The main term is $O(T^{-100})$ if $x<T^{2-3\alpha}$.
\end{lemma}

\begin{lemma}\label{kbeslemma} Let $0<\alpha<\frac{1}{100}$. For any $0<x\le T^{3}$ and any smooth even  function $h$, compactly supported on $(T^{-\alpha}, T^{\alpha})\cup (-T^{\alpha}, -T^{-\alpha})$ with derivatives satisfying $\| h^{(k)} \|_\infty \ll (T^{\alpha})^{k}$, we have
\begin{multline}
\label{kbeslem}  \intt  \sinh(\pi t) K_{2it}(2\pi x) h\Big(\frac{t}{T}\Big) t dt =  \frac{\pi T}{2} \hbar\Big( \frac{\pi x}{T} \Big) - \frac{i\pi^3}{12T}\hbar^{(3)}\Big( \frac{\pi x}{T} \Big)\\+ O\Big( \frac{x}{T^{4-14\alpha}} \Big) +O\Big( \frac{x^2}{T^{5-16\alpha}} \Big) +O(T^{-100}),
\end{multline}
where $\hbar(y)=yh(y)$.
\end{lemma}
\noindent Note that in \cite[Lemma 3.8]{butkha}, the error term $O\big( \frac{x^2}{T^{5-16\alpha}} \big)$ has been erroneously left out. Such a term should be present, as it is for the average of real Bessel functions (see the remark following \cite[Lemma 5.8]{iwaniec}). For our purposes we will have $0<x\le T^{2+4\alpha}$, so the total error for Lemma \ref{kbeslemma} will be essentially the same as for Lemma \ref{jbeslemma}.
 
 \subsection{Kuznetsov's formula for sums of Kloosterman sums}
 
 Let $\Phi$ be a smooth function compactly supported on the positive reals.

 Let 
 \begin{align*}
& \dot{\Phi}(k) = i^k \int_0^\infty J_{k-1}(w)\Phi(w) \frac{dw}{w},\\
& \hat{\Phi}(t)= \frac{i}{2\sinh (\pi t)} \int_0^\infty (J_{2it}(w) - J_{-2it}(w) )\Phi(w) \frac{dw}{w},\\
&\check{\Phi}(t) = \frac{2}{\pi} \cosh(\pi t) \int_0^\infty  K_{2it}(w) \Phi(w) \frac{dw}{w}.
 \end{align*}
 For $q\ge 1$ and $k\ge 2$, let $\mathcal{B}_k(q)$ denote an orthonormal basis of weight $k$ holomorphic cusp forms for $\Gamma_0(q)$. For each element $g$ in this basis, let $n^\frac{k}{2}\rho_g(n)$ denote the $n$-th Fourier coefficient of $g$. Let $\mathcal{B}_0(q)$ denote an orthonormal basis of Maass cusp forms for $\Gamma_0(q)$. For each element $g$ in this basis, let $\rho_g(n)$ denote the $n$-th Fourier coefficient of $g$, and let $\frac{1}{4}+t_g^2$ denote its Laplacian eigenvalue. By a result of Kim and Sarnak \cite[Appendix 2]{kimsar} towards Selberg's Eigenvalue Conjecture, we have that
 \begin{align}
 \label{eigconj} t_g\in \mathbb{R} \ \cup \ \Big(-\frac{7}{64}i, \frac{7}{64}i \Big).
 \end{align}
 Let $\tau_\mathfrak{a}(n,t)$ denote the $n$-th Fourier coefficient of the Eisenstein series $E_\mathfrak{a}(s,\half+it)$ at the cusp $\mathfrak{a}$ of $\Gamma_0(q)$. 
\begin{lemma} \label{kuzback} 
Keep the notation above and let $q\ge 1$. For positive integers $n$ and $m$, we have
\begin{align*}
\sum_{c\ge 1} \frac{S(n,m,qc)}{qc} \Phi\Big( \frac{4\pi\sqrt{nm}}{c}\Big)  = &\sum_{\substack{k\ge 2\\ g \in \mathcal{B}_k(q)}} \dot{\Phi}(k) \frac{(k-1)!\sqrt{nm}}{\pi(4\pi)^{k-1}} \rho_g(n)\overline{\rho_g(m)}\\
+&\sum_{\substack{g \in \mathcal{B}_0(q)}} \hat{\Phi}(t_g) \frac{4\pi\sqrt{nm}}{\cosh(\pi t_g)} \rho_g(n)\overline{\rho_g(-m)}\\
+&\sum_{\mathfrak{a}} \int_{-\infty}^\infty \hat{\Phi}(t)  \frac{\sqrt{nm}}{\cosh(\pi t)} \tau_\mathfrak{a}(n,t)\overline{\tau_\mathfrak{a}(-m,t)} dt,
\end{align*}
and
\begin{align*}
\sum_{c\ge 1} \frac{S(n,-m,qc)}{qc} \Phi\Big( \frac{4\pi\sqrt{nm}}{c}\Big) = 
&\sum_{\substack{g \in \mathcal{B}_0(q)}} \check{\Phi}(t_g) \frac{4\pi\sqrt{nm}}{\cosh(\pi t_g)} \rho_g(n)\overline{\rho_g(m)}\\
+&\sum_{\mathfrak{a}} \int_{-\infty}^\infty \check{\Phi}(t)  \frac{\sqrt{nm}}{\cosh(\pi t)} \tau_\mathfrak{a}(n,t)\overline{\tau_\mathfrak{a}(-m,t)} dt.
\end{align*}
\end{lemma}
\proof
See \cite[Theorem 16.5]{iwakow}, but note that we normalize differently.
\endproof

We now record some properties of the transforms of $\Phi$ given above, based on the situation we will be in (see (\ref{supppsi})). Essentially the same result may be found elsewhere; see for example \cite[Lemma 1]{bloharmic}.
\begin{lemma}\label{psitran}
Suppose that $\Phi(w)$ is supported on 
 \begin{align}
\label{psisupp} X^{-1} T^{-\epsilon} <w<X^{-1}T^{\epsilon}
 \end{align}
 for some $T^{-\epsilon}< X < T^{10}$
  and satisfies
\begin{align}
\label{psicond1}\Phi^{(k)}(w) \ll YX^k T^\epsilon,
 \end{align}
for some $T^{-10}<Y < T^{10}$.

Case 1. If $t\in \mathbb{R}$ and $k\ge 2$ then
 \begin{align}
 \label{psibound} |\dot{\Phi}(k)|, |\hat{\Phi}(t)|, |\check{\Phi}(t)| \ll Y T^\epsilon.
\end{align}
Unless
\begin{align*}
 k,  |t|< T^\epsilon,
\end{align*}
we have
 \begin{align*}
|\dot{\Phi}(k)|, |\hat{\Phi}(t)|, |\check{\Phi}(t)| \ll T^{-B}
\end{align*}
for any $B>0$.
\end{lemma}

Case 2: If $t\in (-\frac{7}{64}i, \frac{7}{64}i )$, then 
 \begin{align}
 \label{psibound2}  |\hat{\Phi}(t)|, |\check{\Phi}(t)| \ll Y X^{\frac{7}{32}} T^\epsilon.
\end{align}

\proof

We demonstrate the claims for $\hat{\Phi}$ only, the other cases being similar. Suppose first that $t \in \mathbb{R} $. Note the bound
\begin{align}
\label{jbound} \frac{J_{2it}(w)}{\sinh(\pi t)}\ll \min\{1,w^{-\half} \},
\end{align}
which follows for $0<w<1$ from the power series \cite[8.402]{gr}
\begin{align}
\label{jpow} J_{\nu}(w)= \sum_{n\ge 0} \frac{(-1)^n (w/2)^{2n+\nu}}{n!\Gamma(n+\nu+1) }
\end{align}
and for $w\ge 1$ by \cite[lemma 6]{blomil}. By (\ref{psisupp}), we may restrict the integral in the definition of $\hat{\Phi}$ to $T^{-11}< w < T^{\epsilon/3}$ and then apply (\ref{psicond1}) with $k=0$ and (\ref{jbound}) to get 
\begin{align*}
|\hat{\Phi}(t)|\ll Y T^\epsilon \int_{T^{-11}}^{T^{\epsilon/3}} \frac{dw}{w}  \ll YT^{\epsilon}.
\end{align*}
This proves (\ref{psibound}). Now suppose that $|t|>T^\epsilon$. By the power series (\ref{jpow}), we have 
\begin{align*}
\hat{\Phi}(t)= \frac{i}{2\sinh (\pi t)} \int_{T^{-11}}^{T^{\epsilon/3}} \Phi(w) \Big( \sum_{n\ge 0} \frac{(-1)^n (w/2)^{2n+2it}}{n!\Gamma(n+2it+1)}  - \sum_{n\ge 0} \frac{(-1)^n (w/2)^{2n-2it}}{n!\Gamma(n-2it+1)}   \Big) \frac{dw}{w} .
\end{align*}
Integrating by parts $k$ times, we get
\begin{align}
\label{serie} \hat{\Phi}(t)&= \frac{i}{2\sinh (\pi t)}   \int_{T^{-11}}^{T^{\epsilon/3}}  \frac{d^k}{dw^k} \Big(\frac{\Phi(w)}{w}\Big) \sum_{n\ge 0} \frac{(-1)^n (w/2)^{2n+2it}w^k}{n!\Gamma(n+2it+1) \prod_{j=1}^{k}(2n+2it+j)}  dw\\
\nonumber & - \frac{i}{2\sinh (\pi t)}   \int_{T^{-11}}^{T^{\epsilon/3}} \frac{d^k}{dw^k} \Big(\frac{\Phi(w)}{w}\Big) \sum_{n\ge 0} \frac{(-1)^n (w/2)^{2n-2it}w^k}{n!\Gamma(n-2it+1) \prod_{j=1}^{k}(2n-2it+j)} dw .
\end{align}
By Stirling's approximation, we have for $|t|>T^\epsilon$ and $0<w<T^{\epsilon/3}$ that
\begin{align*}
\frac{w^{2n}}{\sinh(\pi t)\Gamma(n+1\pm 2it)}\ll \frac{w^{2n}}{|n+1\pm 2it|^{n+\half}}\ll 1.
\end{align*}
By (\ref{psisupp}) we have that
 \begin{align*}
 w^k \frac{d^k}{dw^k} \Big(\frac{\Phi(w)}{w}\Big) \ll  Y X T^\epsilon
 \end{align*}
 for any $k\ge 0$.
Using these bounds in (\ref{serie}) and taking $k$ large, we see that $ |\hat{\Phi}(t)| \ll |t|^{-B}$ for any $B>0$ unless $|t|< T^{\epsilon}$.

Now suppose $t\in (-\frac{7}{64}i, \frac{7}{64}i )$. By \cite[lemma 6]{blomil} and (\ref{jpow}) we have that 
\begin{align*}
 \frac{J_{2it}(w)}{\sinh(\pi t)}\ll \min\{T^\epsilon X^\frac{7}{32},w^{-\half} \} 
\end{align*}
for $w$ in the interval (\ref{psisupp}).
This gives (\ref{psibound2}) by the same argument as above for (\ref{psibound}). 
 \endproof
 
\subsection{Orthonormal basis of newforms}

The right hand side of Lemma \ref{kuzback} involves sums over orthonormal bases of cusp forms. We will need these basis elements to be linear combinations of lifts of newforms. Let $S_k(q)$ denote for $k\ge 2$ the space of holomorphic cusp forms of weight $k\ge 2$ and level $q$, and for $k=0$ the space of Maass cusp forms of level $q$. For $d|q$, let $\mathcal{B}^{*}_k(d)$ denote a basis for the space of newforms of $S_k(d)$, which is orthonormal with respect to the Petersson inner product on $S_k(q)$. For $h\in S_k(d)$ and $b|\frac{q}{d}$, let
\begin{align}
\label{slasdef} h|_b(z) = b^{\frac{k}{2}} h(bz).
\end{align}
\begin{lemma}\label{bases} There exists an orthonormal basis for $S_k(q)$ of the form
\begin{align*}
\bigcup_{d|q} \ \bigcup_{h\in  \mathcal{B}^{*}_k(d)} \{ h_c: c|\tfrac{q}{d} \},
\end{align*}
where
\begin{align*}\
h_c=\sum_{b|c} \kappa_f(c,b) h|_b
\end{align*}
for some numbers $\kappa_f(c,b)\ll q^\epsilon$.
\end{lemma}
\proof See \cite[Lemma 9 and equation (5.6)]{blomil}. This builds on \cite[section 2]{ils} and \cite{rou}.
\endproof

\subsection{Voronoi summation}

The GL(3) Voronoi summation formula was proven by Miller and Schmid \cite{milsch}. Later, Goldfeld and Li \cite{golli} gave another proof and we follow their presentation.

\begin{lemma} \label{vorlem} Let $\psi$ be a smooth, compactly supported function on the positive real numbers. Let $(b,c)=1$ and let $\overline{b}$ denote the multiplicative inverse of $b$ modulo $c$. We have
\begin{multline}
\label{vorstate} \sum_{m\ge 1} A_f(m,r)e\Big(\frac{m\overline{b}}{c}\Big) \psi\Big(\frac{mr^2}{M}\Big)\\ =\sum_\pm \frac{c}{2} \sum_{\substack{\kk \ge 1\\ l|cr }} \frac{A_f(\kk ,l)}{\kk l}S\Big(rb, \pm\kk , \frac{cr}{l}\Big)\frac{1}{2\pi i}\int_{(\sigma)} \Big(\frac{ \kk Ml^2}{c^3r^3}\Big)^{1-s} G^\pm(s)\tilde{\psi}(1-s)ds,
\end{multline}
where $\sigma>0$, $\tilde{\psi}$ denotes the Mellin transform of $\psi$ and 
\begin{align*}
G^\pm(s)=  \frac{\Gam(s+2iT)\Gam(s)\Gam(s-2iT)}{\Gam(1-s-2iT)\Gam(1-s)\Gam(1-s+2iT)}\mp i \frac{\Gam(1+s+2iT)\Gam(1+s)\Gam(1+s-2iT)}{\Gam(2-s-2iT)\Gam(2-s)\Gam(2-s+2iT)}.
\end{align*}
\end{lemma} 
\noindent Suppose that $\|\psi^{(k)} \|_\infty \ll (T^\epsilon)^k$. By integration by parts we can see that $\tilde{\psi}(1-s)\ll (T^{\epsilon })^N (1+|s|)^{-N}$ for any $N\ge 1$. By Stirling's approximation, we have 
\begin{align*}
|G^\pm(s)|\ll (|s+2iT||s||s-2iT|)^{\sigma-\half}.
\end{align*}
Thus we may restrict the integral in (\ref{vorstate}) to $|\Im(s)|<T^\epsilon$. In this range, taking $\sigma=\epsilon$, we get the bound
\begin{align}
\label{hbound} |G^\pm(s)|\ll T^{-1+\epsilon}.
\end{align}
Moving the line of integration in (\ref{vorstate}) far to the right, we see that the sum on the right hand side of (\ref{vorstate}) may be restricted to
\begin{align*}
\kk < \frac{c^3 r^3 T^{2+\epsilon}}{Ml^2},
\end{align*}
up to an error of $O(T^{-100})$. We also observe that in the range $|\Im(s)|<T^\epsilon$, writing $s=\sigma+i\gamma$, Stirling's approximation gives
\begin{align}
\label{vorstirl}& \frac{\Gam(s+2iT)\Gam(s)\Gam(s-2iT)}{\Gam(1-s-2iT)\Gam(1-s)\Gam(1-s+2iT)} = \Big( \frac{T}{\pi}\Big)^{2s-1} \frac{\Gam(s)}{\Gam(1-s)}\Big( 1 + \sum_{n\le N}  \frac{B_n(\sigma,\gamma)}{T^{n}} +O(T^{-\frac{N}{2}})\Big),\\
\nonumber &\frac{\Gam(1+s+2iT)\Gam(1+s)\Gam(1+s-2iT)}{\Gam(2-s-2iT)\Gam(2-s)\Gam(2-s+2iT)} =  \Big( \frac{T}{\pi}\Big)^{2s-2} \frac{\Gam(1+s)}{\Gam(2-s)}\Big( 1 + \sum_{n\le N}  \frac{C_n(\sigma,\gamma)}{T^{n}} +O(T^{-\frac{N}{2}})\Big)
\end{align} 
for any $N\ge 1$ and some $B_n(\sigma,\gamma)$ and $C_n(\sigma,\gamma)$ polynomial in $\sigma$ and $\gamma$.
\section{Proof of Lemma \ref{eislem}.}
It is more convenient to renormalize $f$ so that $\|f\|_2=1$. This does affect what needs to be proved, for in Lemma \ref{contpart} we need only an upper bound. So let
\begin{align}
\label{foucusp} f(x+iy)=\rho_f(1)\sum_{n\neq 0} \lambda_f(n)\sqrt{y}K_{iT}(2\pi n y)e(nx)
\end{align}
denote the Fourier series expansion of $f$, where
\begin{align*}
\rho_f(1)^2 = \frac{2\cosh(\pi T)}{L(1,\sym^2 f)}.
\end{align*}
By unfolding, we have for $\Re(s)>1$,
\begin{align*}
\langle  f^2, E(\cdot,s) \rangle = 2\rho_f(1)^2 \sum_{n\ge 1} \frac{\lambda_f(n)^2}{(2\pi n)^{s}} \int_0^\infty  y^{s} (K_{iT}(y))^2 \frac{dy}{y}.
\end{align*}
We have that
\begin{align*}
\sum_{n\ge 1} \frac{\lambda_f(n)^2}{(2\pi n)^{s}} = \frac{\zeta(s)L(s,\sym^2 f)}{(2\pi)^s \zeta(2s)}
\end{align*}
by \cite[Chapter 7]{gol},
\begin{align*}
\int_0^\infty y^s (K_{iT}(y))^2 \frac{dy}{y} = \frac{2^{s-3} \Gamma^2(\tfrac{s}{2})\Gamma(\tfrac{s}{2}+iT)\Gamma(\tfrac{s}{2}-iT)}{\Gamma(s)}
\end{align*}
by \cite[6.576]{gr}, and
\begin{align*}
\cosh(\pi T)=\frac{\pi}{|\Gamma(\half+iT)|^2}.
\end{align*}
Thus taking $s=\half+it$ by analytic continuation, we have
\begin{align}
\label{unfold} |\langle f^2, E(\cdot,\thalf +it) \rangle|^2 =  \frac{\pi}{4} H(t) \frac{|\zeta(\half+it) L(\thalf+it,\sym^2 f)|^2}{L(1,\sym^2 f)^2 |\zeta(1+2it)|^2}.
\end{align}
This may be compared with (\ref{wats}).

By (\ref{stirht}), we may restrict the integral in Lemma \ref{eislem} to $|t|<T^{1+\epsilon}$. The zeta and $L$-value in the denominator of (\ref{unfold}) are on the edge of the region of absolute convergence, so they are bounded below by $T^{-\epsilon}$. Thus 
\begin{align*}
 \int_{-\infty}^{\infty} |\langle f^2,  E(\cdot,\thalf +it) \rangle|^2 dt \ll T^{\epsilon} \int_{-T^{1+\epsilon}}^{T^{1+\epsilon}} \frac{|\zeta(\thalf+it) L(\thalf+it,\sym^2 f)|^2}{(1+|t|)\prod_\pm(1+|2T\pm t|)^\half} dt.
\end{align*}
On the GLH, this is bounded by
\begin{align*}
T^{\epsilon} \int_{-T^{1+\epsilon}}^{T^{1+\epsilon}} \frac{1}{(1+|t|)\prod_\pm(1+|2T\pm t|)^\half} dt \ll T^{-1+\epsilon}.
\end{align*}

\section{Proof of Proposition \ref{prop}: Applying the trace formula}

We first refine what needs to be proved for Proposition \ref{prop}. We can immediately treat the contribution to (\ref{remains}) of $t_j$ close to zero and close to $2T$. Let $0<\alpha<\frac{1}{100}$ be a fixed constant to be determined later. On the GLH, we have by (\ref{wats}) and (\ref{stirht}) that
\begin{align*}
\sum_{|t_j|<T^{1-\alpha}}  + \sum_{|t_j-2T|<T^{1-\alpha}} | \langle f^2 , u_j \rangle |^2 \ll T^{-1+\epsilon} \sum_{|t_j|<T^{1-\alpha}}  \frac{1}{1+|t|} + T^{-\frac{3}{2}+\epsilon}\sum_{|t_j-2T|<T^{1-\alpha}} \frac{1}{(1+|t_j-2T|)^\half}.
\end{align*} 
By Weyl's law (see \cite[page 391]{iwakow}) we have that this is less than $T^{-\frac{\alpha}{2}+\epsilon}$. Thus we may restrict the left hand side of (\ref{remains}) to values of $t_j$ that are roughly of size $2T$, but not too close to $2T$. We have to take care when making this restriction because we must use functions that will satisfy the conditions of Kuznetsov's trace formula. 

\begin{lemma} \label{wlem} Let $0<\alpha<\frac{1}{100}$ be a fixed constant to be determined later and define the even function
\begin{align*}
W(t) =W_\alpha(t) = \left(1-\exp \left(-\left(\frac{t}{(2T)^{1-\frac{\alpha}{2}}}\right)^{2\lceil \frac{1000}{\alpha} \rceil} \right) \right) \left(1-\exp \left(-\left(\frac{4T^2-t^2}{4T^{2-\frac{\alpha}{2}}}\right)^{2\lceil \frac{1000}{\alpha} \rceil} \right) \right),
\end{align*}
where $\lceil x \rceil$ denotes the least integer greater than or equal to $x$.
We have that $H(t)W(t)\ll T^{-100}$ unless 
\begin{align}
\label{wrange} T^{1-\alpha} <|t| <2T - T^{1-\alpha},
\end{align}
in which range
\begin{align}
\label{hhbound}  \frac{d^k}{dt^k} H(t)W(t) \ll T^{-2} \left( T^{-1+\alpha} \right)^k
\end{align}
for $k \ge 0$. In the range $T^{1-\frac{\alpha}{4}} < |t| <2T - T^{1-\frac{\alpha}{4}}$, we have that $W(t)=1+O(T^{-100})$. 
\end{lemma}
\proof
Suppose that $|t|\le T^{1-\alpha}$. Then
\begin{align*}
 \left(\frac{t}{(2T)^{1-\frac{\alpha}{2}}}\right)^{2\lceil \frac{1000}{\alpha} \rceil} \ll T^{-100}\implies 1-\exp \left(-\left(\frac{t}{(2T)^{1-\frac{\alpha}{2}}}\right)^{2\lceil \frac{1000}{\alpha} \rceil} \right) \ll T^{-100}
\end{align*}
and trivially,
\begin{align*}
1-\exp \left(-\left(\frac{4T^2-t^2}{4T^{2-\frac{\alpha}{2}}}\right)^{2\lceil \frac{1000}{\alpha} \rceil} \right) \ll 1.
\end{align*}
Therefore $H(t)W(t)\ll W(t)\ll T^{-100}$ for $|t|\le T^{1-\alpha}$. Now suppose $|t|\ge 2T - T^{1-\alpha}$. If also $|t|>2T+T^{\epsilon}$, then by (\ref{stirht}) we have $H(t)\ll T^{-100}$, and so $W(t)H(t)\ll T^{-100}$. So suppose  $2T-T^{1-\alpha}\le |t|\le 2T+T^\epsilon$. Then 
\begin{align*}
\left(\frac{4T^2-t^2}{4T^{2-\frac{\alpha}{2}}}\right)^{2\lceil \frac{1000}{\alpha} \rceil} \ll T^{-100 }\implies 1-\exp \left(-\left(\frac{4T^2-t^2}{4T^{2-\frac{\alpha}{2}}}\right)^{2\lceil \frac{1000}{\alpha} \rceil}\right)\ll T^{-100}
\end{align*}
and trivially,
\begin{align*}
 1-\exp \left(-\left(\frac{t}{(2T)^{1-\frac{\alpha}{2}}}\right)^{2\lceil \frac{1000}{\alpha} \rceil} \right) \ll 1. 
\end{align*}
This proves the first claim. 

For the second claim, observe that for $x>0$ and $N>0$, we have $\frac{d^k}{dx^k} \exp(-x^N) \ll 1$. So 
\begin{align}
\label{wder} W^{(k)}(t) \ll \left( T^{-1+\alpha} \right)^k.
\end{align}
To prove the same sort of bound for $H(t)$, we need more terms in the Stirling expansion (\ref{stirht}). In the range (\ref{wrange}), we have
\begin{align}
\label{h0def} H(t)=\frac{1}{2T^2} H_0\Big(\frac{|t|}{2T}\Big),
\end{align}
where
\begin{multline*}
H_0(x) = \frac{8\pi }{ x(1-x^2)^{\half}}\\\times  \Big(1+\sum_{n\le N} \frac{B_n}{T^n}\Big)\Big(1+\sum_{n\le N} \frac{C_n}{(Tx)^n}\Big)\Big(1+\sum_{n\le N}\frac{D_n}{(2T(1+x))^n}\Big)\Big(1+\sum_{n\le N}\frac{D_n}{(2T(1-x))^n}\Big)+O( T^{-\frac{N}{2}})
\end{multline*}
for some constants $B_n, C_n, D_n$. By taking $N$ large enough we see that if $T^{-\alpha}\ll |x| \ll 1-T^{-\alpha}$, then
\begin{align}
\label{h0der}H_0^{(k)}(x) \ll (T^\alpha)^k .
\end{align}
Thus
\begin{align}
H^{(k)}(t) \ll T^{-2} \left( T^{-1+\alpha} \right)^k
\end{align}
in the range (\ref{wrange}).

For the third claim, suppose that $T^{1-\frac{\alpha}{4}} < |t| <2T - T^{1-\frac{\alpha}{4}}$.  Then
\begin{align*}
 \left(\frac{t}{(2T)^{1-\frac{\alpha}{2}}}\right)^{2\lceil \frac{1000}{\alpha} \rceil} \gg T^{100}\implies 1-\exp \left(-\left(\frac{t}{(2T)^{1-\frac{\alpha}{2}}}\right)^{2\lceil \frac{1000}{\alpha} \rceil} \right) =1 + O(T^{-100})
\end{align*}
and
\begin{align*}
\left(\frac{4T^2-t^2}{4T^{2-\frac{\alpha}{2}}}\right)^{2\lceil \frac{1000}{\alpha} \rceil} \gg T^{100 }\implies 1-\exp \left(-\left(\frac{4T^2-t^2}{4T^{2-\frac{\alpha}{2}}}\right)^{2\lceil \frac{1000}{\alpha} \rceil}\right)= 1+O(T^{-100}).
\end{align*}

\endproof

Proposition \ref{prop} is thus reduced to proving
\begin{align}
\label{prop2} \sum_{j\ge 1} | \langle  f^2 , u_j \rangle|^2 W(t_j) = \frac{2\pi }{3}+O(T^{-\delta})
\end{align}
for some $\delta>0$.

We remark that since $L(\half, u_j)$ and $L(\half, u_j\times \sym^2 f)$ are non-negative (by the work of \cite{lapral} and \cite{lap}), the contribution of the very small eigenvalues is at least
\begin{align*}
\sum_{|t_j|<T^{\epsilon}}   | \langle f^2 , u_j \rangle |^2 \gg T^{-\epsilon} \sum_{|t_j|<T^{\epsilon}}  \frac{L(\half, u_j)L(\half,u_j\times \sym^2 f)}{T}.
\end{align*} 
Thus even if we were not assuming GLH, a subconvexity bound for $L(\half, u_j\times \sym^2 f)$ in the $T$ aspect, which is polynomial in $|t_j|<T^{\epsilon}$, would be required, but this is an unsolved and very difficult problem.

By (\ref{wats}), we have
\begin{align}
\label{watsapplied} \sum_{j\ge 1} | \langle f^2 , u_j \rangle |^2 W(t_j) =  \frac{\pi^3}{72 L(1,\sym^2 f)^2} \sum_{j \ge 1} H(t_j)W(t_j)\frac{L(\half,u_j) L(\half, u_j \times \sym^2 f )}{L(1,\sym^2 u_j)}.
\end{align}
For $u_j$ even, we may use the approximate functional equation (\ref{afe2}) to write
\begin{align*}
L(\thalf,u_j) L(\thalf, u_j \times \sym^2 f ) = 2 L(\thalf, u_j) \sum_{m,r\ge 1} \frac{\lambda_j(m)A_f(r,m)}{rm^\half}V_2(r^2m,t_j). 
\end{align*}
But this equality holds for $u_j$ odd as well, since in this case both sides vanish. Now we may use the approximate functional equation (\ref{trickafe}) for $L(\half,u_j)$, which holds for both even and odd forms. This idea is an important feature of our proof which, as mentioned in the remarks following Lemma \ref{afelem}, will make the analysis more pleasant. Thus we get that the right hand side of (\ref{watsapplied}) equals
\begin{align}
 \label{pleas} &\sum_{\pm} \frac{\pi^3}{72 L(1,\sym^2 f)^2} \sum_{j\ge 1}  \sum_{n\ge 1} \frac{\lambda_j(\pm n)}{n^{\half}} L(u_j\times \sym^2 f) H(t_j)W(t_j)V_1^\pm (n,t_j)V_2(r^2m,t_j)\\
\nonumber &+O\Big(T^{-\half+\frac{\beta}{2}+\alpha+\epsilon} \sum_{T^{1-\alpha}<t_j<T^{1+\epsilon}} |H(t_j) W(t_j)| \Big| \sum_{m,r\ge 1} \frac{\lambda_j(m)A_f(r,m)}{rm^\half}V_2(r^2m,t_j) \Big|\Big)
\end{align}
By (\ref{hhbound}) and (\ref{videf}), the error term is
\begin{align*}
O\Big(T^{-\frac{5}{2}+\frac{\beta}{2}+\alpha+\epsilon} \sum_{T^{1-\alpha}<t_j<T^{1+\epsilon}} \Big|  \int_{(\epsilon)} L(\thalf+s, u_j\times \sym^2 f) e^{s^2}  \frac{G_2(\half+s)}{G_2(\half)} \frac{ds}{s}  \Big|\Big),
\end{align*}
which is $O(T^{-\half+\frac{\beta}{2}+\frac{\alpha}{2}+\epsilon})$ on the GLH. By the approximate functional equation (\ref{afe2}), the main term of (\ref{pleas}) equals
\begin{align*}
\sum_{\pm} \frac{\pi^3}{36 L(1,\sym^2 f)^2} \sum_{n,m,r\ge 1} \frac{A_f(m,r)}{r(nm)^{\half}}\sum_{j\ge 1} \lambda_j(\pm n)\lambda_j(m) H(t_j)W(t_j)V_1^\pm (n,t_j)V_2(r^2m,t_j).
\end{align*}
Applying the Kuznetsov trace formula to the inner sum gives:
\begin{align}
\sum_{j\ge 1} | \langle f^2 , u_j \rangle |^2 W(t_j) =  \mathcal{D} + \mathcal{E} + \mathcal{O}^+ + \mathcal{O}^- +O(T^{-\delta})
\end{align}
for some $\delta>0$, where
\begin{align}
&\label{ddef} \mathcal{D} = \frac{\pi^3}{36L(1,\sym^2 f)^2}  \sum_{n,r\ge 1} \frac{A_f(n,r)}{rn}\int_{-\infty}^{\infty} H(t)W(t)V_1^\pm (n,t)V_2(r^2n,t) \frac{d^*t}{2\pi^2},\\
&\nonumber \mathcal{E}=   \frac{-\pi^3}{18 L(1,\sym^2 f)^2} \int_{-\infty}^{\infty} \sum_\pm \sum_{n,m,r\ge 1} \frac{A_f(m,r)\lambda(n,t)\lambda(m,-t)}{r(nm)^{\half} |\zeta(1+2it)|^2}  H(t)W(t)V_1^\pm (n,t)V_2(r^2m,t) \frac{dt}{2\pi}
\end{align}
and 
\begin{multline*}
\mathcal{O}^\pm = \frac{\pi^3}{36L(1,\sym^2 f)^2}   \sum_{n,m,r\ge 1} \sum_{c\ge 1} \frac{A_f(m,r)}{r(nm)^{\half}} \frac{S(\pm n, m, c)}{c}\\  \times \int_{-\infty}^{\infty} \mathcal{J}^{\pm}\Big(\frac{\sqrt{nm}}{c},t\Big) H(t)W(t)V_1^\pm (n,t)V_2(r^2m,t) \frac{d^*t}{2\pi}.
\end{multline*}
By the decay of $V_1^\pm$ and $V_2$ , we may restrict the sum above to $n<T^{1\mp \beta+\epsilon}$ and $mr^2<T^{3+\epsilon}$. We may also restrict to $c\le T^3$, say, by a standard method (see \cite[Lemma 5]{blo} for example). To see this, consider $\mathcal{O}^+$ for instance, and move the line of the $t$-integral from $\Im(t)=0$ to $\Im(t)=-\frac{1}{2}+\epsilon$. On the new line, we have by the power series (\ref{jpow}) that 
\begin{align*}
\mathcal{J}^{+}\Big(\frac{\sqrt{nm}}{c},t\Big)\ll \Big(\frac{\sqrt{nm}}{c}\Big)^{1-\epsilon}
\end{align*}
 for $\frac{\sqrt{nm}}{c}<1$, as is the case when $c> T^3$. The rest of the integrand satisfies
\begin{align*}
H(t)W(t)V_1^\pm (n,t)V_2(r^2m,t)\tanh(\pi t)\ll \exp(-\pi q(t,T)),
\end{align*}
where $q(t,T)$ is given in (\ref{qtT}). Thus the contribution of the terms with $c>T^3$ is bounded by
\begin{align*}
\sum_{\substack{n\le T^{1+\epsilon}\\mr^2\le T^{3+\epsilon} }} \sum_{c>T^3} \frac{|A_f(m,r)|}{r(nm)^{\half}}   \frac{|S(\pm n, m, c)|}{c} \int_{-T^{1+\epsilon}}^{T^{1+\epsilon}} \Big(\frac{\sqrt{nm}}{c}\Big)^{1-\epsilon} dt \ll T^{1+\epsilon}  \sum_{c>T^3}  \frac{1}{c^{\frac{3}{2}-\epsilon}}. 
\end{align*}
The last bound uses the average version of the Ramanujan bound given in (\ref{avgram}) and Weil's bound for the Kloosterman sum. The result is less than a negative power of $T$.

We will prove that $\mathcal{D}$ yields the main term, while $\mathcal{E}$ and $\mathcal{O}^\pm$ are bounded by a negative power of $T$.

\section{Proof of Proposition \ref{prop}:  The diagonal}

The goal of this section (see (\ref{prop2}) is to show that
\begin{align*}
\mathcal{D} = \frac{2\pi }{3}+O(T^{-\delta}),
\end{align*}
for some $\delta>0$. 
This was sketched in \cite[Section 4]{blokhayou} but here we provide the details. By (\ref{stirv1}-\ref{stirv2}) and (\ref{ddef}), we have
\begin{multline}
\label{dexp} \mathcal{D} =  \frac{\pi^3}{36L(1,\sym^2 f)^2} \int_{-\infty}^{\infty} H(t)W(t) \frac{1}{(2\pi i)^2} \int_{(\epsilon)} \int_{(\epsilon)} e^{s_1^2+s_2^2} T^{-\beta s_1} \left( \frac{|t|}{2\pi}\right)^{s_1}  \left( \frac{|t(4T^2-t^2)|}{8\pi^3}\right)^{s_2}\\
\times \sum_{n,r\ge 1} \frac{A_f(n,r)}{r^{1+2s_2}n^{1+s_1+s_2}} \frac{ds_1}{s_1} \frac{ds_2}{s_2} \frac{d^*t}{2\pi^2}\Big(1 + O(T^{-1+\alpha})\Big),
\end{multline}
where $0<\alpha<\frac{1}{100}$ is as in Lemma \ref{wlem}. 
By \cite[Proposition 6.6.3]{gol}, we have 
\begin{align*}
\sum_{n,r\ge 1} \frac{A_f(n,r)}{r^{1+2s_2}n^{1+s_1+s_2}} = \frac{L(1+2s_2,\sym^2 f)L(1+s_1+s_2,\sym^2 f)}{\zeta(2+s_1+3s_2)}.
\end{align*}
Thus assuming the GLH, we have that the error term in (\ref{dexp}) contributes
\begin{align*}
O\left(T^{-1+\alpha+\epsilon} \int_{-\infty}^{\infty} |H(t)W(t)| t \tanh(\pi t) dt \right).
\end{align*}
From Lemma \ref{wlem}, the weight function $H(t)W(t)$ is $O(T^{-100})$ unless $T^{1-\alpha} < |t| <2T - T^{1-\alpha}$. By the estimate for $H(t)$ given in \ref{stirht}), we have
\begin{align*}
\int_{T^{1-\alpha} }^{2T - T^{1-\alpha}} |H(t)W(t) | t \tanh (\pi t) dt \ll  \int_{0 }^{2T } \frac{1}{(4T^2-t^2)^\half} dt \ll 1
\end{align*}
and so the error term of (\ref{dexp}) is $O(T^{-\half} )$.

Consider the main term of (\ref{dexp}). Moving the line of integration to $\Re(s_1)=-\frac{1}{10}$, we pick up a simple pole at $s_1=0$, getting
\begin{multline}
\label{1stshift} \mathcal{D} =  \frac{\pi^3}{36L(1,\sym^2 f)^2} \int_{-\infty}^{\infty} H(t)W(t) \frac{1}{(2\pi i)^2} \int_{(\epsilon)}  e^{s_2^2}   \left( \frac{|t(4T^2-t^2)|}{8\pi^3}\right)^{s_2}\\
\times \frac{L(1+2s_2,\sym^2 f)L(1+s_2,\sym^2 f)}{\zeta(2+3s_2)} \frac{ds_2}{s_2} \frac{d^*t}{2\pi^2}  +  O\left(T^{-\frac{1}{10}(1-\beta)+\epsilon}\right).
\end{multline}
The new error term arises by applying GLH on the shifted line of integration, and it is $O(T^{-\frac{1}{20}})$. Now moving the line of integration to $\Re(s_2)=-\frac{1}{10}$, and picking up a simple pole at $s_2=0$, we get 
\begin{align*}
\mathcal{D} =  \frac{\pi^3}{36 L(1,\sym^2 f)^2} \int_{-\infty}^{\infty} \frac{L(1,\sym^2 f)^2}{\zeta(2)} H(t)W(t) \frac{d^*t}{2\pi^2}  + O(T^{-\frac{1}{10}(3-\alpha)+\epsilon} + T^{-\frac{1}{20}} ).
\end{align*}
The error term is $O(T^{-\frac{1}{20}})$. The main term equals
\begin{align*}
 \frac{1}{6 \pi } \int_{0}^{\infty} H(t)W(t) t\tanh(\pi t) dt.
\end{align*}
We can now restrict the integrand to the range $T^{1-\alpha} < t <2T - T^{1-\alpha}$, on which interval $\tanh(\pi t) = 1+O(T^{-100})$.  Further, for $T^{1-\frac{\alpha}{4}} < t < 2T- T^{1-\frac{\alpha}{4}}$ we have $W(t) = 1+O(T^{-100})$. Thus
\begin{align*}
\mathcal{D} &= \frac{1}{6 \pi } \int_{T^{1-\alpha}}^{2T-T^{1-\alpha}}  H(t)W(t) t dt +O(T^{-\frac{1}{20}})\\
&= \frac{1}{6 \pi } \int_{T^{1-\frac{\alpha}{4}}}^{2T-T^{1-\frac{\alpha}{4}}}  H(t)t dt +O(T^{-\frac{\alpha}{8}}).
\end{align*}
By (\ref{stirht}) we have
\begin{align*}
\mathcal{D}&= \frac{1}{6 \pi } \int_{T^{1-\frac{\alpha}{4}}}^{2T-T^{1-\frac{\alpha}{4}}}  \frac{8\pi}{(4T^2-t^2)^\half} dt + O(T^{-\frac{\alpha}{8}}) \\
&= \frac{4}{3} \left( \arcsin\Big( \frac{2T-T^{1-\frac{\alpha}{4}}}{2T} \Big)-\arcsin\Big( \frac{T^{1-\frac{\alpha}{4}}}{2T} \Big) \right) + O(T^{-\frac{\alpha}{8}}))\\
&= \frac{4}{3}\Big(\frac{\pi}{2}-0\Big) + O(T^{-\frac{\alpha}{8}}).
\end{align*}
The error term is some negative power of $T$. This completes the evaluation of the diagonal.

\section{Proof of Proposition \ref{prop}:  The Eisenstein series contribution}

In this section we show that $\mathcal{E}$ is bounded by a negative power of $T$. We first rewrite the expression for $\mathcal{E}$ using the following approximate functional equations, which are analogous to those in Lemma \ref{afelem}:
\begin{align*}
&|\zeta(\thalf + it)|^2= \sum_\pm \sum_{n\ge 1} \frac{\lambda(n,t)}{n^\half} V_1^\pm(n,t),\\
&|L(\thalf + it,\sym^2 f)|^2=   \sum_{m,r\ge 1} \frac{A_f(m,r)\lambda(m,-t)}{rm^\half} V_2(r^2m,t).
\end{align*}
We get that
\begin{align*}
\mathcal{E}=   \frac{-\pi^3}{9 L(1,\sym^2 f)^2} \int_{0}^{\infty} \frac{ |\zeta(\thalf + it)|^2 |L(\thalf + it,\sym^2 f)|^2}{  |\zeta(1+2it)|^2}  H(t)W(t)  \frac{dt}{2\pi}.
\end{align*}
By the decay of the weight function $H(t)W(t)$, we may restrict the integral to $T^{1-\alpha} < t <2T - T^{1-\alpha}$, and then apply the GLH and the estimate for $H(t)$ given in (\ref{stirht}) to see that 
\begin{align*}
\mathcal{E} &=   \frac{-\pi^3}{9 L(1,\sym^2 f)^2} \int_{T^{1-\alpha}}^{ 2T - T^{1-\alpha} } \frac{ |\zeta(\thalf + it)|^2 |L(\thalf + it,\sym^2 f)|^2}{  |\zeta(1+2it)|^2}  H(t)W(t)  \frac{dt}{2\pi}\\
&\ll T^{\epsilon} \int_{T^{1-\alpha}}^{ 2T - T^{1-\alpha} } \frac{1}{t(4T^2-t^2)^{\half}} dt\\
&\ll T^{-1+\alpha+\epsilon}.
\end{align*}

\section{Proof of Proposition \ref{prop}:  The short off-diagonal} 

The goal of this section is to show that $\mathcal{O}^+$ is bounded by a negative power of $T$. We have seen that for any $0<\alpha<\frac{1}{100}$, we have $H(t)W(t)\ll T^{-100}$ unless $T^{1-\alpha}<|t|<2T -T^{1-\alpha}$. Thus in the expression for $\mathcal{O}^+$ we may restrict the integral to this range. We may also replace $d^*t$ by $tdt$ because $\tanh t= 1 +O(T^{-100})$ in the given range of $t$. Thus
\begin{multline*}
\mathcal{O}^+ \ll T^\epsilon  \sum_{\substack{ n <T^{1-\beta+\epsilon} \\ mr^2<T^{3+\epsilon} }}  \sum_{c\le T^3} \frac{A_f(m,r)}{r(nm)^{\half}} \frac{S(n, m, c)}{c}\\
\times  \int_{T^{1-\alpha}}^{2T-T^{1-\alpha} } + \int_{-2T+T^{1-\alpha}}^{-T^{1-\alpha} }  \mathcal{J}^{+}\Big(\frac{\sqrt{nm}}{c},t\Big) H(t)W(t)V_1^+ (n,t)V_2(r^2m,t) tdt.
\end{multline*}

Let $Z$ be any smooth, even function compactly supported on $(T^{-\alpha},2 - T^{-\alpha})\cup (-2 + T^{-\alpha},-T^{-\alpha})$ with derivatives satisfying
\begin{align*}
\| Z^{(k)} \|_\infty \ll (T^\alpha)^k.
\end{align*}
Then $Z(\frac{t}{T})$ is supported on $T^{1-\alpha}<|t|<2T -T^{1-\alpha}$, on which $\frac{d^k}{dt^k} Z(\frac{t}{T})\ll T^{(-1+\alpha)k}$, and we may use such a function to approximate the characteristic function of this interval. We may absorb $W(t)$ into the function $Z(\frac{t}{T})$ by property (\ref{wder}). Writing 
\begin{align*}
H(t)=\frac{1}{2T^2}H_0\Big(\frac{t}{T}\Big)
\end{align*}
as in (\ref{h0def}), we may also absorb $H_0(\frac{t}{T})$ into $Z(\frac{t}{T})$ by property (\ref{h0der}). Thus it suffices to prove that for any function $Z$ as above, we have
\begin{align*}
\frac{1}{T^{2}} \sum_{n,m,r\ge 1}  \sum_{c\le T^3} \frac{A_f(m,r)}{r(nm)^{\half}} \frac{S(n, m, c)}{c} 
  \int_{-\infty}^{\infty }  \mathcal{J}^{+}\Big(\frac{\sqrt{nm}}{c},t\Big) Z\Big(\frac{t}{T}\Big) V_1^+ (n,t)V_2(r^2m,t) tdt \ll T^{-\delta} 
\end{align*}
for some $\delta>0$. We apply Lemma \ref{jbeslemma} to evaluate the Bessel transform. Note that the function $Z(\frac{t}{T})V_1^+ (n,t)V_2(r^2m,t)$ satisfies the conditions of the lemma by the remarks above and by (\ref{vder1}-\ref{vder2}). 

The contribution of the main term of Lemma \ref{jbeslemma} is $O(T^{-100})$, unless
\begin{align*}
\frac{2\sqrt{nm}}{c} > T^{2-3\alpha}. 
\end{align*}
Since by the decay of  $V_1^+ (n,t)V_2(r^2m,t)$ we may take $n <T^{1-\beta+\epsilon}$ and $m<T^{3+\epsilon}$ up to $O(T^{-100})$, this imposes
\begin{align*}
c <  T^{3\alpha-\frac{\beta}{2}}.
\end{align*}
We now fix
\begin{align*}
\beta=7\alpha
\end{align*}
so that the above condition on $c$ is impossible and the main term is $O(T^{-100})$. This of course leads to an acceptable bound for $\mathcal{O}^+$.

The error term $O(\frac{\sqrt{nm}}{cT^{3-12\alpha}})$ arising from Lemma \ref{jbeslemma} contributes
\begin{align*}
O\Big(\frac{1}{T^{2}} \sum_{\substack{ n <T^{1-7\alpha+\epsilon} \\ mr^2<T^{3+\epsilon} }} \sum_{c\le T^3} \frac{|A_f(m,r)|}{r(nm)^{\half}} \frac{|S(n, m, c)|}{c} \frac{\sqrt{nm}}{cT^{3-12\alpha}}\Big)
\end{align*}
By (\ref{avgram}) and Weil's bound for the Kloosterman sum, this is
\begin{align*}
O\Big(T^{-5+12\alpha+\epsilon} \sum_{\substack{ n <T^{1-7\alpha+\epsilon} \\ mr^2<T^{3+\epsilon} }} \sum_{c\le T^3}  \frac{1}{rc^{\frac{3}{2}}}  \Big).
\end{align*}
The innermost $c$-sum is $O(1)$, so the line above is $O(T^{-1+5\alpha+ \epsilon})$, which is admissible as we assume $\alpha<\frac{1}{100}$.

\section{Proof of Proposition \ref{prop}:  The long off-diagonal} 

The goal now is to show that $\mathcal{O}^-$ is bounded by a negative power of $T$. This proof is the heart of our paper. As in the previous section, it suffices to prove that for any smooth, even function $Z$ compactly supported on 
\begin{align*}
(T^{-\alpha},2 - T^{-\alpha})\cup (-2 + T^{-\alpha},-T^{-\alpha})
\end{align*}
 with derivatives satisfying $\| Z^{(k)} \|_\infty \ll (T^\alpha)^k$, we have 
\begin{align*}
\frac{1}{T^{2}} \sum_{\substack{ n <T^{1+7\alpha+\epsilon} \\ mr^2<T^{3+\epsilon} }}  \sum_{c\le T^3}  \frac{A_f(m,r)}{r(nm)^{\half}} \frac{S(-n, m, c)}{c} 
  \int_{-\infty}^{\infty }  \mathcal{J}^{-}\Big(\frac{\sqrt{nm}}{c},t\Big) Z\Big(\frac{t}{T}\Big) V_1^- (n,t)V_2(r^2m,t) tdt \ll T^{-\delta} 
\end{align*}
for some $\delta>0$. We may replace $V_1^- (n,t)$ and $V_2(r^2m,t)$ by the main terms in their Stirling expansions (\ref{stirv1}) and (\ref{stirv2}) since the lower order terms can be treated similarly. Thus we need to show that 
\begin{align}
\label{suff} \frac{1}{T^{2}} \sum_{\substack{ n <T^{1+7\alpha+\epsilon} \\ mr^2<T^{3+\epsilon} }}  \sum_{c\le T^3}  \frac{A_f(m,r)}{r(nm)^{\half}} \frac{S(-n, m, c)}{c} 
  \int_{-\infty}^{\infty }  \mathcal{J}^{-}\Big(\frac{\sqrt{nm}}{c},t\Big) Z\Big(\frac{t}{T}\Big) V\Big(\frac{n}{T},\frac{mr^2}{T^3},\frac{t}{T}\Big) tdt 
\end{align}
is bounded by a negative power of $T$ where
\begin{align*}
V\left(x_1,x_2; y \right) =  \int_{(\sigma)} e^{s_1^2} ( 2\pi x_1)^{-s_1} |T^{7\alpha}y|^{s_1}   \frac{ds_1}{s_1} \cdot \int_{(\sigma)} e^{s_2^2} (8\pi^3 x_2  )^{-s_2} |y(4-y^2)|^{-s_2}   \frac{ds_2}{s_2}
\end{align*}
for any $\sigma>0$.
We apply Lemma \ref{kbeslemma} to evaluate the Bessel transform in (\ref{suff}). The error term arising from this result contributes less than a negative power of $T$, just as in the previous section. There are two similar, non-oscillatory main terms in the asymptotic given by Lemma \ref{kbeslemma}. It suffices to treat only the leading main term as the other will contribute a factor of $T$ less. Therefore the goal is to bound by a negative power of $T$ the sum
\begin{align*}
\frac{1}{T^{2}} \sum_{\substack{ n <T^{1+7\alpha+\epsilon} \\ mr^2<T^{3+\epsilon} }}  \sum_{c\le T^3}  \frac{A_f(m,r)}{r} \frac{S(-n, m, c)}{c^2}  Z\Big(\frac{2\pi\sqrt{nm}}{Tc}\Big) V\Big(\frac{n}{T},\frac{mr^2}{T^3},\frac{2\pi\sqrt{nm}}{Tc}\Big).
\end{align*}
Applying a smooth partition of unity, we consider the sum above in dyadic intervals. Let $U$ be a smooth bump function supported on $(1, 2)\times(1, 2)$ and possessing bounded derivatives. It suffices to prove that
\begin{align}
\label{ready} \frac{1}{T^{2}} \sum_{n,m,r\ge 1}  \sum_{c\ge 1} \frac{A_f(m,r)}{r} \frac{S(-n, m, c)}{c^2}  Z\Big(\frac{2\pi\sqrt{nm}}{Tc}\Big) V\Big(\frac{n}{T},\frac{mr^2}{T^3},\frac{2\pi\sqrt{nm}}{Tc}\Big)U\Big(\frac{n}{N},\frac{mr^2}{M}\Big)\ll T^{-\delta}
\end{align}
for some $\delta>0$, where
\begin{align*}
N<T^{1+7\alpha+\epsilon}, M<T^{3+\epsilon},  \frac{\sqrt{NM}}{rT} \ll c\ll \frac{\sqrt{NM}}{rT^{1-\alpha}}.
\end{align*}
The bounds on $c$ are enforced by the function $Z$, and make the condition $c\le T^3$ redundant. It will be apparent from the proof that as long as $\alpha$ is smaller than some fixed constant, there exists some absolute $\delta>0$ independent of the value of $\alpha$. For this reason it will be very convenient to rename $\alpha$ to $\epsilon$ and apply the $\epsilon$-convention. Therefore we have
\begin{align}
\label{ranges} N<T^{1+\epsilon}, M<T^{3+\epsilon},  \frac{\sqrt{NM}}{rT} \ll c\ll \frac{\sqrt{NM}}{rT^{1-\epsilon}}.
\end{align}

\subsection{Poisson and Voronoi summation}

By separating the $n$-sum in (\ref{ready}) into residue classes mod $c$ and applying the Poisson summation formula,  we get that the left hand side of (\ref{ready}) equals
\begin{multline}
\label{pois} \frac{N}{T^{2}}  \sum_{\substack{m,r,c\ge 1}}  \sum_{a\bmod c}     \frac{A_f(m,r)}{r} \frac{S(-a, m, c)}{c^3} \\ \times \sum_{-\infty<\jj  <\infty} e\Big(\frac{\jj  a}{c}\Big) \int_{-\infty}^{\infty} U\Big(\xi ,\frac{mr^2}{M}\Big) Z\Big(\frac{2\pi\sqrt{\xi Nm}}{Tc}\Big) V\Big(\frac{\xi N}{T},\frac{mr^2}{T^3},\frac{2\pi\sqrt{\xi N m}}{Tc}\Big) e\Big( \frac{-\jj  N \xi }{c} \Big) d\xi.
\end{multline}
Now in this expression, writing $x=\frac{mr^2}{M}$, we may replace $m$ by $\frac{Mx}{r^2}$. Writing $y=\frac{\sqrt{NM}}{crT}$, we may replace $N$ by $\frac{y^2c^2r^2T^2}{M}$. Thus (\ref{pois}) equals
\begin{align}
\label{pois2} \frac{N}{T^{2}}  \sum_{\substack{c,m,r\ge 1\\-\infty<\jj  <\infty}}  \sum_{a\bmod c} \frac{A_f(m,r)}{r} \frac{S(-a, m, c)}{c^3}  e\Big(\frac{\jj  a}{c}\Big)  \psi\Big(\frac{mr^2}{M},\frac{\sqrt{NM}}{crT};  \frac{\jj  cr^2T^2  }{M} \Big) 
\end{align}
where
\begin{align}
\label{psidef} \psi(x,y; u) =  \int_{-\infty}^{\infty} U(\xi,x )  Z(2\pi y\sqrt{\xi x } ) V\Big(\frac{\xi N}{T},\frac{xM}{T^3}, 2\pi y \sqrt{\xi x} \Big) e(-\xi  y^2 u) d\xi.
\end{align}
In the integrand above, $x$ and $\xi$ must lie in the interval $(1,2)$ by the definition of $U$, while $y$ must lie in the interval $(T^{-\epsilon},T^{\epsilon})$ by the definition of $Z$. Thus all three variables $x,\xi,y$ should be thought of as roughly constant, and bounded away from zero.

We first observe that we may restrict (\ref{pois2}) to $\jj \neq 0$. First note that the contribution of $\jj=0$ is nil unless $c=1$, because $\sum_{a\bmod c} S(-a,m,c)=0$ for $c>1$. This leaves the case $\jj=0$ and $c=1$, whose contribution is
\begin{align}
\nonumber&\frac{N}{T^{2}}  \sum_{\substack{m,r\ge 1}}  \frac{A_f(m,r)}{r}  \psi\Big(\frac{mr^2}{M},\frac{\sqrt{NM}}{rT}; 0 \Big)\\
\nonumber&=  \sum_{\substack{m,r\ge 1}}  \frac{A_f(m,r)}{mr}  \frac{mr^2}{M} \frac{NM}{r^2T^2} \psi\Big(\frac{mr^2}{M},\frac{\sqrt{NM}}{rT}; 0 \Big)\\
\label{j0c1} &= \sum_{\substack{m,r,v\ge 1}} \frac{\mu(v)}{v^2} \frac{A_f(m,1)A_f(r,1)}{mr}  \frac{mr^2v^3}{M} \frac{NM}{r^2 v^2 T^2} \psi\Big(\frac{mr^2 v^3}{M},\frac{\sqrt{NM}}{r vT}; 0 \Big),
\end{align}
where the last equality follows by the Hecke relations (\ref{gl3hecke}). The sum is trivially $O(T^{\epsilon})$, using (\ref{avgram}), so we must save any negative power of $T$. We are done unless
\begin{align*}
v\ll T^{\epsilon} \ \ \text{ and } \ \  \frac{mr^2v^3}{M}\gg T^{-\epsilon},
\end{align*}
in which case, since  
\begin{align*}
\frac{mr^2v^3}{M}\ll T^\epsilon \ \ \text{ and } \ \  T^{-\epsilon}\ll  \frac{NM}{r^2 v^2 T^2}  \ll T^{\epsilon}
\end{align*}
are enforced by the weight function,
we have
\begin{align*}
\frac{T^{2-\epsilon}}{N} \ll  m \ll \frac{T^{2+\epsilon}}{N}
\end{align*}
This shows that the interval of summation of $m$ is at least as long as $T^{1-\epsilon}$, by (\ref{ranges}). On the GLH we can show in a standard way (see section \ref{glhsection}) that
\begin{align*}
\sum_{m\le x} A_f(m,1) \ll x^{\frac{1}{2}+\epsilon}.
\end{align*}
Using this and partial summation gives the required saving in (\ref{j0c1}).

The next step is to transform the $m$-sum using GL(3) Voronoi summation (see Lemma \ref{vorlem}). 
Writing
\begin{align*}
S(-a, m, c) = \summ_{b \bmod c} e\Big( \frac{-ab+m\overline{b}}{c} \Big), 
\end{align*}
we get by Voronoi summation that the part of (\ref{pois2}) with $\jj\neq 0$ equals
\begin{multline}
\label{vor} \frac{N}{T^2} \sum_{\pm} \sum_{\substack{c,r,\kk ,|\jj  |\ge 1\\ l|cr}}  \frac{A_f(\kk ,l)}{ \kk l r c^2} \summ_{b \bmod c}  S\Big(rb, \pm\kk , \frac{cr}{l}\Big) \sum_{a \bmod c} e\Big( \frac{a(\jj-b)}{c} \Big)   \\
\times \frac{1}{4\pi i} \int_{(\epsilon)}  \Big(\frac{ \kk Ml^2}{c^3r^3}\Big)^{1-s} G^\pm(s)   \int_0^\infty  \psi\Big(x,\frac{\sqrt{NM}}{crT};  \frac{\jj  cr^2T^2  }{M} \Big)    x^{-s}  dx \ ds,
\end{multline}
where $G^\pm(s)$ was defined in Lemma \ref{vorlem}. We have that $\sum_{a \bmod c} e( \frac{a(\jj-b)}{c} )$ equals $c$ if $b \equiv \jj  \bmod c$, and 0 otherwise. In the case $b \equiv \jj  \bmod c$, we have $br \equiv \jj  r \bmod cr$. Thus (\ref{vor}) equals
\begin{multline}
\label{vor2} \frac{N}{T^2}  \sum_{\pm} \sum_{\substack{c,r,\kk ,|\jj  |\ge 1\\ (\jj ,c)=1\\ l | cr}}  \frac{A_f(\kk ,l)}{\kk lrc}   S\Big(\jj  r, \pm\kk , \frac{cr}{l}\Big) \\ \times \frac{1}{4\pi i} \int_{(\epsilon)}  \Big(\frac{ \kk Ml^2}{c^3r^3}\Big)^{1-s} G^\pm(s)  \int_0^\infty   \psi\Big(x,\frac{\sqrt{NM}}{crT};  \frac{\jj  cr^2T^2  }{M} \Big)    x^{-s}  dx \ ds.
\end{multline}
Recall by the remarks following Lemma \ref{vorlem}, that up to negligible error, we can restrict the sum to
\begin{align}
\label{kcap} \frac{ \kk Ml^2}{c^3r^3}<T^{2+\epsilon}.
\end{align}
So if we set $z=\frac{ \kk Ml^2}{c^3r^3T^2}$ then $T^{-6} <z<T^{\epsilon}$ and we can write $M= \frac{zc^3r^3T^2} {\kk l^2}$ in the fraction $ \frac{\jj  cr^2T^2  }{M}$ appearing in (\ref{vor2}), and $\kk =\frac{zc^3r^3T^2}{ Ml^2}$. Thus (\ref{vor2}) equals
\begin{align}
\label{vor3} \frac{NM}{T^2} \sum_{\pm} \sum_{\substack{c,r,\kk ,|\jj  |\ge 1\\ (\jj ,c)=1\\ l | cr}} \frac{l A_f(\kk ,l) }{ c^4r^4 }   S\Big(r \jj , \pm\kk , \frac{cr}{l}\Big)  \Psi^\pm \Big(\frac{\sqrt{NM}}{crT} ,\frac{ \kk Ml^2}{c^3r^3T^2} ;\frac{\jj \kk l^2}{c^2 r} \Big),
\end{align}
where
\begin{align}
\label{bigpsidef} \Psi^\pm (y,z;u) = \frac{1}{4\pi i} \int_{(\epsilon)} z^{-s} T^{-2s} G^\pm(s)   \int_0^\infty   \psi (x,y;  z^{-1} u )    x^{-s}  dx  \ ds.
\end{align}
The goal is to show that (\ref{vor3}) is bounded by a negative power of $T$. Using $G^\pm(s)\ll T^{-1+\epsilon}$, a bound given in (\ref{hbound}), and repeatedly integrating by parts in (\ref{psidef}), we have
\begin{align}
\label{psider2} \frac{\partial^k}{\partial u^k} \Psi^\pm (y,z; u) \ll T^{-1+\epsilon} z^{-k} |z^{-1} u|^{-A}
\end{align}
for any $k\ge 0$ and $A\ge 0$. Thus up to negligible error, we may assume
\begin{align}
\label{cap} |u| < zT^\epsilon.
\end{align}
Using this, we have
\begin{align}
\label{psider1} \frac{\partial^k}{\partial y^k} \Psi^\pm (y,z;u) \ll T^{-1+\epsilon}, \ \ \frac{\partial^k}{\partial z^k} \Psi^\pm (y,z;u) \ll T^{-1+\epsilon}z^{-k}
\end{align}

\subsection{Preparation for Kuznetsov's formula}

We first explain the idea of what we are about to do next. Our task is to bound by a negative power of $T$ a sum like (when $r=l=1$),
\begin{align}
\label{sumidea} \frac{NM}{T^2} \sum_{\substack{c,\kk,|\jj  |\ge 1 }}  \frac{A_f(\kk ,1)}{c^4}   S\Big(\jj , \pm\kk , c\Big)  \Psi^\pm \Big(\frac{\sqrt{NM}}{cT} ,\frac{\kk M}{c^3T^2} ;\frac{\jj \kk }{c^2 } \Big).
\end{align}
By (\ref{cap}) and $z<T^\epsilon$, we see that up to negligible error, we may assume $|\jj  |\kk \ll c^2T^\epsilon$. Also recall that $| \Psi^\pm (y,z;u)|\ll T^{-1+\epsilon}$. By these remarks, the trivial bound for (\ref{sumidea}) is $O(T^{\half+\epsilon})$, and now we must save this much and a little more. We have already exploited the sums over $\jj$ and $\kk$ (they arose through Poisson and Voronoi summation). Now we will exploit the sum over $c$. Since the range of $c$ in the sum is at least as large as $T^{-\epsilon} \sqrt{|\jj  |\kk }$, a range sometimes referred to as the Linnik range (see \cite{sartsi}), we are in a good position to use Kuznetsov's formula to transform the sum of Kloosterman sums 
 into a sum of Hecke eigenvalues. This final sum will be estimated under the GLH to complete the proof. To get a feel of how this works, consider the generic ranges $c\asymp T, N\asymp T, M\asymp T^3$, that imply $\jj  \ll T^\epsilon, \kk \asymp T^2$. Using Kuznetsov's formula, we will get an identity of the shape
 \begin{align*}
 \sum_{c\asymp \sqrt{|\jj |\kk} } c^{-1} S(\pm \jj , \kk, c) \asymp \sum_{t_j <T^{\epsilon}}\lambda_j(\jj )\lambda_j(\kk)+\ldots
 \end{align*}
 Thus (\ref{sumidea}) essentially becomes, ignoring the $\jj$ sum of length $T^\epsilon$ and keeping in mind that $\Psi^\pm(y,z,w)$ has a factor $T^{-1}$ by (\ref{vorstirl}), 
  \begin{align*}
\sum_{t_j <T^{\epsilon}} T^{-2} \sum_{\kk\asymp T^2} A_f(\kk,1)\lambda_j(\kk)+\ldots
 \end{align*}
 On the GLH, we get cancellation in the $\kk$-sum and obtain the required bound.

To carry out the above program, we must first take care of the technicalities posed by the presence of the $l$ and $r$ parameters. First we detect the condition $(\jj ,c)=1$ in (\ref{vor3}) using the M\"{o}bius function  by recalling that $\sum_{\substack{d|c, d|\jj}} \mu(d)$ equals 1 if $(\jj ,c)=1$ and 0 otherwise. Thus we need to bound the following sum, for each sign $\pm$, by a negative power of $T$:
\begin{align*}
\frac{NM}{T^2} \sum_{\substack{d,c,r,\kk ,|\jj  |\ge 1 \\ l | cdr}}  \frac{l\mu(d) A_f(\kk ,l)}{c^4r^4d^4}   S\Big(r\jj  d, \pm\kk , \frac{cdr}{l}\Big)  \Psi^\pm \Big(\frac{\sqrt{NM}}{cdrT} ,\frac{ \kk Ml^2}{c^3d^3r^3T^2} ;\frac{\jj \kk l^2}{c^2 d r} \Big).
\end{align*}
We can reorder this sum by the gcd of $l$ and $c$ to say that it equals
\begin{align*}
\frac{NM}{T^2} \sum_{b\ge 1} \  \sum_{\substack{d,c,r,\kk ,|\jj  |\ge 1 \\ l | cdr \\ (l,c)=b}} \frac{l\mu(d) A_f(\kk ,l)}{c^4r^4d^4}    S\Big(r\jj  d, \pm\kk , \frac{cdr}{l}\Big)  \Psi^\pm \Big(\frac{\sqrt{NM}}{cdrT} ,\frac{ \kk Ml^2}{c^3d^3r^3T^2} ;\frac{\jj \kk l^2}{c^2 d r} \Big).
\end{align*}
Replacing $l$ by $lb$ and $c$ by $cb$, this sum equals
\begin{align*}
\frac{NM}{T^2} \sum_{b\ge 1} \  \sum_{\substack{d,c,r,\kk ,|\jj  |\ge 1 \\ l | dr \\ (l,c)=1}}  \frac{l\mu(d) A_f(\kk ,bl)}{c^4r^4d^4 b^3}   S\Big(r\jj  d, \pm\kk , \frac{cdr}{l}\Big)  \Psi^\pm \Big(\frac{\sqrt{NM}}{cbdrT} ,\frac{ \kk Ml^2}{c^3d^3r^3bT^2} ;\frac{\jj \kk l^2}{c^2 d r} \Big).
\end{align*}
Now detecting $(l,c)=1$ using the M\"{o}bius function, the sum above equals
\begin{align*}
\frac{NM}{T^2}  \sum_{\substack{b,d,c,r,\kk ,|\jj  |\ge 1 \\ al | dr }}  \frac{l\mu(d)\mu(a) A_f(\kk ,abl)}{c^4r^4d^4 a^3b^3}   S\Big(r\jj  d, \pm\kk , \frac{cdr}{l}\Big)  \Psi^\pm \Big(\frac{\sqrt{NM}}{cabdrT} ,\frac{ \kk Ml^2}{c^3d^3r^3abT^2} ;\frac{\jj \kk l^2}{c^2 d r} \Big).
\end{align*}
Using the Hecke relations (\ref{gl3hecke}) and writing $\frac{1}{(cabdr)^3}=(\frac{\sqrt{NM}}{cabdrT})^3(\frac{T}{\sqrt{NM}})^3$, it suffices, to show that
\begin{multline}
\label{kuzprep} \frac{T}{\sqrt{NM}}  \sum_{\substack{b,d,r \le T^{2} \\ al | dr \\ v|abl}} |A_f(abl, 1)|  \Bigg|  \sum_{\substack{\kk,|\jj|,c\ge 1}} \frac{A_f(\kk ,1)}{\frac{cdr}{l}} S\Big(dr \jj , \pm\kk  v, \frac{cdr}{l}\Big) \\ \Big(\frac{\sqrt{NM}}{cabdrT}\Big)^3  \Psi^\pm \Big(\frac{\sqrt{NM}}{cabdrT} ,\frac{ \kk Ml^2}{c^3d^3r^3ab T^2} ;\frac{\jj \kk  v l^2}{c^2 d r } \Big) \Bigg| \ll T^{-\delta}
\end{multline}
for some $\delta>0$. The given bounds $b,d,r \le T^{2}$ are enforced by the weight functions. Also note for later that the desired bound of (\ref{kuzprep}) is trivial for the sub-sum restricted to $|\jj \kk| \le T^{\frac{1}{10}}$ using $\Psi^\pm(y,z,u)\ll T^{-1+\epsilon}$ and $\frac{T}{\sqrt{NM}}\ll  \frac{T^{\epsilon}}{cabdr}$, and Kim and Sarnak's bound
\begin{align*}
A_f(abl,1)\ll (abl)^{\frac{7}{32}+\epsilon}.
\end{align*}

We continue to reshape the inner sum in order to apply Lemma \ref{kuzback}. We write
\begin{multline*}
 \sum_{c\ge 1} \frac{S\Big(dr\jj , \pm\kk  v, \frac{cdr}{l}\Big)}{\frac{cdr}{l}}   \Big(\frac{\sqrt{NM}}{cabdrT}\Big)^3 \Psi^\pm \Big(\frac{\sqrt{NM}}{cabdrT} ,\frac{ \kk Ml^2}{c^3d^3r^3abT^2} ; \frac{\jj \kk  v l^2}{c^2 d r} \Big) \\= \sum_{c\ge 1} \frac{S\Big(r\jj  d, \pm\kk  v, \frac{cdr}{l}\Big)}{\frac{cdr}{l}} \Phi\Big(\frac{4\pi\sqrt{|\jj  |\kk  rdv}}{\frac{cdr}{l}}\Big),
\end{multline*}
where, after absorbing the factor $(\frac{\sqrt{NM}}{cabdrT} )^3$ into the function $\Psi^\pm$, we let
\begin{align*}
\Phi(w)=   \Psi^\pm \Big( w \frac{\sqrt{NM}}{4\pi abl T   \sqrt{|\jj  |\kk rdv}} , \ w^3 \frac{ \kk M }{ab l T^2(4\pi \sqrt{|\jj  |\kk rdv})^3} ; \ w^2   \frac{\text{sgn}(j)}{(4\pi)^2}\Big)
\end{align*}
for $w>0$. Since $\Psi^{\pm}(y,z;u)$ is restricted to $T^{-\epsilon}\ll y \ll T^{\epsilon}$, we have
\begin{align}
\label{supppsi} \text{support}(\Phi(w)) \subset \Big( T^{-\epsilon}\Big(\frac{\sqrt{NM}}{ abl T   \sqrt{|\jj  |\kk rdv}}\Big)^{-1} ,  \Big(\frac{\sqrt{NM}}{ abl T   \sqrt{|\jj  |\kk rdv}}\Big)^{-1}T^{\epsilon} \Big).
\end{align}
Further, we have seen that $\Psi^{\pm}(y,z;u)$ may be restricted to $z \ll T^{\epsilon}$ and $u\ll zT^{\epsilon}$ up to an error of $O(T^{-100})$. Thus by (\ref{psider2}) and (\ref{psider1}), we have
\begin{align*}
 |\Phi^{(k)}(w)|\ll T^{-1+\epsilon} (1+w^{-k} )  (z^{-1}w^2)^{-A} + O(T^{-100})
\end{align*}
for any $k\ge0$ and $A\ge0$. Equivalently, by (\ref{supppsi}) we have
\begin{align}
\label{phik2} |\Phi^{(k)}(w)|\ll \frac{1}{T^{1-\epsilon}} \Big(1+\Big( \frac{\sqrt{NM}}{ abl T   \sqrt{|\jj  |\kk rdv}}\Big) ^{k} \Big)\Big(\frac{z}{w^2}\Big)^{A}+ O(T^{-100}).
\end{align}
We see that up to an error of $O(T^{-100})$ say, we can assume
\begin{align*}
w^2\ll zT^{\epsilon} \ \Rightarrow \ \frac{l\sqrt{|\jj  |\kk  rdv}}{cdr} \ll T^{\epsilon}.
\end{align*}
So, since $\frac{\sqrt{NM}}{cabdrT}$ is restricted to $(T^{-\epsilon},T^{\epsilon})$, we have
\begin{align}
\label{xlower} \frac{\sqrt{NM}}{ abl T   \sqrt{|\jj  |\kk rdv}}\gg T^{-\epsilon}.
\end{align}
Thus (\ref{phik2}) implies 
\begin{align}
\label{phik3} |\Phi^{(k)}(w)|\ll \frac{1}{T^{1-\epsilon}} \Big( \frac{\sqrt{NM}}{ abl T   \sqrt{|\jj  |\kk rdv}}\Big) ^{k} + O(T^{-100}).
\end{align}

Finally we note some bounds on $\jj$ and $\kk$. By (\ref{cap}), we may restrict to
\begin{align*}
\frac{\jj \kk v l^2}{c^2 d r} \ll \frac{\kk M l^2}{c^3 d^3 r^3 ab T^{2-\epsilon}} , 
\end{align*}
which gives
\begin{align*}
\jj \ll \frac{M}{cd^2 r^2abv T^{2-\epsilon}}.
\end{align*}
We already have the bound (\ref{kcap}) on $\kk$. So since $T^{-\epsilon}< \frac{\sqrt{NM}}{cabdrT} < T^{\epsilon}$, we may eliminate $c$ and record the bounds
\begin{align*}
&|\jj| \ll  \frac{ M^\half }{T^{1-\epsilon} N^\half drv},\\
& \kk \ll \frac{N^\frac{3}{2} M^\half}{ a^2b^2 l^2 T^{1-\epsilon}}.
\end{align*}

Also note for later use that $\Phi(w)$ depends implicitly on $\jj$ and $\kk$, and by (\ref{psider1}) we have
\begin{align}
\label{jkders} \frac{\partial^n}{\partial \jj^n} \Phi(w) \ll T^{-1+\epsilon}\jj^{-n}, \ \ \  \frac{\partial^n}{\partial \kk^n} \Phi(w) \ll T^{-1+\epsilon}\kk^{-n}.
\end{align}

\subsection{Kuznetsov's trace formula and the large sieve}

Consider the sub-sum of (\ref{kuzprep}) consisting of the terms with $\jj>0$ and the positive sign case, the rest of the sum being similarly treated. Applying Lemma \ref{kuzback} and the remarks following it, we see that it suffices to show that the following expression is bounded by a negative power of $T$:
\begin{align}
\label{exp} &\frac{T}{\sqrt{NM}} \sum_{\substack{b,d,r \le T^{2} \\ al | dr \\ v|abl}}  (abl)^{\frac{7}{32}} \Bigg|  \sum_{\substack{K\le \kk < 2K \\ J\le \jj< 2J}}  A_f(\kk ,1)  \sum_{\substack{g \in \mathcal{B}_0(\frac{dr}{l})\\ |t_g|<T^\epsilon}} \hat{\Phi}(t_g) \frac{4\pi\sqrt{\jj\kk drv}}{\cosh(\pi t_g)} \rho_g(\kk v)\overline{\rho_g(\jj dr)} \Bigg| \\
\nonumber &+ \frac{T}{\sqrt{NM}} \sum_{\substack{b,d,r \le T^{2} \\ al | dr \\ v|abl}} (abl)^{\frac{7}{32}} \Bigg|  \sum_{\substack{K\le \kk < 2K \\ J\le \jj< 2J}} A_f(\kk ,1) \sum_{\substack{g \in \mathcal{B}_k(\frac{dr}{l})\\ k<T^\epsilon }} \dot{\Phi}(k) \frac{(k-1)!\sqrt{\jj\kk drv}}{\pi(4\pi)^{k-1}}\rho_g(\kk v)\overline{\rho_g(\jj dr)}\Bigg| \\
\nonumber &+\frac{T}{\sqrt{NM}} \sum_{\substack{b,d,r \le T^{2} \\ al | dr \\ v|abl}}  (abl)^{\frac{7}{32}} \Bigg|  \sum_{\substack{K\le \kk < 2K \\ J\le \jj< 2J}} A_f(\kk ,1) \sum_{\mathfrak{a}} \int_{|t|<T^\epsilon} \hat{\Phi}(t)  \frac{\sqrt{\jj\kk drv}}{\cosh(\pi t)} \tau_\mathfrak{a}(\kk v,t)\overline{\tau_\mathfrak{a}(-\jj dr,t)} dt\Bigg|
\end{align}
for any positive integers
\begin{align}
\label{jk} J< \frac{ M^\half }{T^{1-\epsilon} N^\half drv}, \ \ \ K< \frac{N^\frac{3}{2} M^\half}{ a^2b^2 l^2 T^{1-\epsilon}}
\end{align}
with 
\begin{align*}
JK>T^\frac{1}{10}.
\end{align*}
For this last assumption, see the remark following (\ref{kuzprep}).
We only treat the first line of (\ref{exp}) as the rest are similar.

Let
\begin{align*}
 X = \frac{\sqrt{NM}}{ T^{1-\epsilon}abl   \sqrt{JK rdv}},
\end{align*}
so that $\Phi(w)$ is supported on $X^{-1}\ll w\ll X^{-1}$, by (\ref{supppsi}).
We first reduce to the case 
\begin{align*}
bdr \ll T^{\epsilon} \ \text{ and } \ X \ll T^{\epsilon}.
\end{align*}
{\bf Non-exceptional eigenvalues.} Consider the contribution of $t_g\in \mathbb{R}$ to the first line of (\ref{exp}). We have the bound
\begin{align*}
|\hat{\Phi}(t_g)|\ll T^{-1+\epsilon},
\end{align*}
 by (\ref{phik3}) and (\ref{psibound}). Using the Cauchy-Schwarz inequality, and enlarging the spectral sum to $|t_g|<X$ (recall that $X\gg T^{\epsilon}$ by (\ref{xlower})), we would like to say that the contribution of the non-exceptional eigenvalues to the first line of (\ref{exp}) is bounded by 
 \begin{align}
\label{cs1} \frac{T^\epsilon}{\sqrt{NM}}\sum_{\substack{b,d,r \le T^{2} \\ al | dr \\ v|abl}} (abl)^{\frac{7}{32}}   \frac{l}{dr}   \Big(\sum_{\substack{g \in \mathcal{B}_0(\frac{dr}{l})\\ |t_g|<X}} \Big| \sum_{\jj  } \alpha_\jj   \rho_g(\jj) \Big|^2 \Big)^\half 
  \Big(\sum_{\substack{g \in \mathcal{B}_0(\frac{dr}{l})\\ |t_g|<X}} \Big| \sum_{\kk } \beta_\kk   \rho_g(\kk) \Big|^2 \Big)^\half, 
 \end{align}
 where 
 \begin{align*}
 \alpha_\jj  =
 \begin{cases}
 \Big(\frac{dr}{l} \frac{4\pi \jj  }{ \cosh(\pi t_g )} \Big)^\half &\text{ for } Jdr\le \jj < 2Jdr, \  \jj\equiv 0 \bmod dr, \\
 0 &\text{ otherwise},
 \end{cases}
 \end{align*}
 and 
  \begin{align*}
 \beta_\kk  =
 \begin{cases}
 A(\frac{\kk}{v},1) \Big(\frac{dr}{l} \frac{4\pi \kk }{ \cosh(\pi t_g)} \Big)^\half &\text{ for } Kv\le \kk < 2Kv, \ \kk \equiv 0 \bmod v, \\
 0 &\text{ otherwise}.
 \end{cases}
 \end{align*}
 However, although (\ref{cs1})  is essentially the right upper bound, to make the argument rigorous we must remember that $\hat{\Phi}(t_g)$ depends implicitly on $\jj$ and $\kk$. We can separate variables as follows. Write
 \begin{align*}
 \hat{\Phi}(t_g) = \hat{\Phi}(t_g,\jj,\kk).
 \end{align*}
 Insert the factors $U_1(\frac{\jj}{J})$ and $U_2(\frac{\kk}{K})$, where $U_1(x_1)$ and $U_2(x_2)$ are smooth bump functions that are compactly supported on $(\frac{1}{2},\frac{3}{2})$ and equal 1 on $[ 1,2 ]$. Using the Mellin transform, we have
 \begin{multline*}
 U_1\Big(\frac{\jj}{J}\Big) U_2\Big(\frac{\kk}{K}\Big) \hat{\Phi}(t_g,\jj,\kk) \\
 = \frac{1}{(2\pi i)^2} \int_{(\epsilon)} \int_{(\epsilon)} \Big(\frac{\jj}{J}\Big)^{-s} \Big(\frac{\kk}{K}\Big)^{-s} \int_{\frac{1}{2}}^{\frac{3}{2}} \int_{\frac{1}{2}}^{\frac{3}{2}} x_1^{s_1} x_2^{s_2} U_1(x_1) U_2(x_2) \hat{\Phi}(t_g,x_1J,x_2K) \ \frac{dx_1}{x_1}\frac{dx_2}{x_2} \ ds_1ds_2.
  \end{multline*}
  We can restrict the integrals to $|\Im(s_1)|,|\Im(s_2)|<T^\epsilon$ by integrating by parts (using (\ref{jkders})).
 
 By the spectral large sieve \cite[Theorem 7.24, equation (7.40)]{iwakow} and (\ref{avgram2}), we have that (\ref{cs1}) is bounded by
 \begin{align}
\label{sls} \frac{T^\epsilon}{\sqrt{NM}} \sum_{\substack{b,d,r \le T^{2}  \\ al | dr \\ v|abl}} (abl)^{\frac{7}{32}}  \frac{l}{dr}  \Big( \Big(\frac{dr}{l} X^2 + Jdr\Big)J \Big)^\half \Big( \Big(\frac{dr}{l} X^2 + K v\Big)K \Big)^\half.
 \end{align}
 We expand this out and look at the cross terms one by one. We have 
 \begin{align}
\label{c1}  \frac{T^\epsilon}{\sqrt{NM}} \sum_{\substack{b,d,r \le T^{2} \\ al | dr \\ v|abl}} (abl)^{\frac{7}{32}}  X^2 (JK)^\half \ll \sum_{\substack{b,d,r \le T^{2} \\ al | dr \\ v|abl}} (abl)^{\frac{7}{32}}   \frac{(NM)^\half}{T^{2-\epsilon} a^2b^2l^2rdv(JK)^\half}.
 \end{align}
 Now because $N\ll T^{1+\epsilon}$ and $M\ll T^{3+\epsilon}$ and we assume $JK>T^\frac{1}{10}$, this is bounded by a negative power of $T$. The next cross term is  
  \begin{align}
\label{c2}  \frac{T^\epsilon}{\sqrt{NM}} \sum_{\substack{b,d,r \le T^{2} \\ al | dr \\ v|abl}} (abl)^{\frac{7}{32}}  \Big(\frac{l}{dr}\Big)^\half X J^\half K v^\half \ll \sum_{\substack{b,d,r \le T^{2} \\ al | dr \\ v|abl}} (abl)^{\frac{7}{32}}  \frac{N^\frac{3}{4}M^\frac{1}{4} }{rd a^2 b^2 l^\frac{3}{2} T^{\frac{3}{2}-\epsilon}}.
 \end{align}
Thus the left hand side is bounded by a negative power of $T$ unless $N\gg T^{1-\epsilon}$, $M\gg T^{3-\epsilon}$, $abl\ll T^{\epsilon}$, and
\begin{align*}
K\gg T^{-\epsilon} \frac{N^\frac{3}{2} M^\half}{ a^2b^2 l^2 T}\gg T^{2-\epsilon}.
\end{align*}
 The next cross term is
  \begin{align}
 \label{c3}  \frac{T^\epsilon}{\sqrt{NM}} \sum_{\substack{b,d,r \le T^{2} \\ al | dr \\ v|abl}} (abl)^{\frac{7}{32}}  \Big(\frac{l}{dr}\Big)^\half X K^\half J (dr)^\half \ll \sum_{\substack{b,d,r \le T^{2} \\ al | dr \\ v|abl}} (abl)^{\frac{7}{32}}  \frac{M^\frac{1}{4}}{T^{\frac{3}{2}-\epsilon}N^\frac{1}{4} ab l^\half drv}.
 \end{align}
 This is less than a negative power of $T$ since $M\ll T^{3+\epsilon}$. The final cross term is
   \begin{align}
\label{c4} \frac{T^\epsilon}{\sqrt{NM}} \sum_{\substack{b,d,r \le T^{2} \\ al | dr \\ v|abl}}(abl)^{\frac{7}{32}}  \frac{lv^\half}{(dr)^\half} JK \ll  \sum_{\substack{b,d,r \le T^{2} \\ al | dr \\ v|abl}} (abl)^{\frac{7}{32}}  \frac{(NM)^\half}{T^{2-\epsilon} a^2b^2l (dr)^\frac{3}{2} v^\half}.
 \end{align}
So the left hand side is bounded by a negative power of $T$ unless $N\gg T^{1-\epsilon}$, $M\gg T^{3-\epsilon}$, $abdr\ll T^{\epsilon}$ (which implies $l\ll T^{\epsilon}$), and 
\begin{align*}
K\gg T^{-\epsilon} \frac{N^\frac{3}{2} M^\half}{ a^2b^2 l^2 T}\gg T^{2-\epsilon}.
\end{align*}
 
 We conclude that (\ref{sls}) is bounded by a negative power of $T$ unless 
 \begin{align*}
 T^{1-\epsilon}\ll N \ll T^{1+\epsilon}, \ \ T^{3-\epsilon}\ll M\ll T^{3+\epsilon}, \ \ b\ll T^{\epsilon}, \ \ T^{2-\epsilon}\ll K \ll T^{2+\epsilon},
 \end{align*}
 which is the case we consider now. In these ranges, we deduce from (\ref{jk}) that we must have
 \begin{align*}
J\ll T^{\epsilon}, \ \ dr \ll T^\epsilon.
 \end{align*}
 The last bound holds because were the contrary true, we would have $J\ll T^{-\epsilon}$ and that is impossible for a positive integer. It follows that we also have $a,l,v\ll T^{\epsilon}$, since these are divisors of small quantities. In summary, we can assume
  \begin{align}
\label{summary1} T^{1-\epsilon}\ll N \ll T^{1+\epsilon}, \ \ T^{3-\epsilon}\ll M\ll T^{3+\epsilon}, \ \ J\ll T^{\epsilon}, \ \  T^{2-\epsilon}\ll K \ll T^{2+\epsilon}, \ \ ablvdr\ll T^\epsilon.
 \end{align}
  From this we get
 \begin{align}
 \label{summary2} X\ll T^\epsilon.
 \end{align}
 
 {\bf Exceptional eigenvalues.} We now consider the contribution of $t_g\in (-\frac{7}{64}i, \frac{7}{64}i)$ to the first line of (\ref{exp}). This time have the bound
\begin{align*}
|\hat{\Phi}(t_g)|\ll T^{-1+\epsilon}X^\frac{7}{32}
\end{align*}
 by (\ref{phik3}) and (\ref{psibound2}). By the Cauchy-Schwarz inequality, the contribution of the possible exceptional eigenvalues is (essentially) bounded by 
 \begin{align*}
 \frac{T^\epsilon}{\sqrt{NM}}\sum_{\substack{b,d,r \le T^{2} \\ al | dr \\ v|abl}} (abl X)^{\frac{7}{32}}   \frac{l}{dr}  \Big(\sum_{\substack{g \in \mathcal{B}_0(\frac{dr}{l})\\ |t_g|<1}} \Big| \sum_{\jj  } \alpha_\jj   \rho_g(\jj) \Big|^2 \Big)^\half 
  \Big(\sum_{\substack{g \in \mathcal{B}_0(\frac{dr}{l})\\ |t_g|<1}} \Big| \sum_{\kk } \beta_\kk   \rho_g(\kk) \Big|^2 \Big)^\half, 
 \end{align*}
 where $ \alpha_\jj $ and $ \beta_\kk $ are as above. By the spectral large sieve again, this is bounded by
 \begin{align}
 \label{sls2} \frac{T^\epsilon}{\sqrt{NM}} \sum_{\substack{b,d,r \le T^{2}  \\ al | dr \\ v|abl}} (ablX)^{\frac{7}{32}}  \frac{l}{dr}  \Big( \Big(\frac{dr}{l}  + Jdr\Big)J \Big)^\half \Big( \Big(\frac{dr}{l}  + K v\Big)K \Big)^\half.
 \end{align}
 We know from (\ref{xlower}) that $X>T^{\epsilon}$. Expand (\ref{sls2}) and consider each cross term. The cross terms corresponding to (\ref{c1}-\ref{c3}) are clearly smaller in size, because each of (\ref{c1}-\ref{c3}) has a factor of $X$ or $X^2$, while the cross terms of (\ref{sls2}) have a factor of $X^{\frac{7}{32}}$ only. Thus it remains only to consider the cross term
 \begin{align*}
  \frac{T^\epsilon}{\sqrt{NM}} \sum_{\substack{b,d,r \le T^{2} \\ al | dr \\ v|abl}}(abl X)^{\frac{7}{32}}  \frac{lv^\half}{(dr)^\half} JK.
   \end{align*}
   Since $X>1$, this is less than
   \begin{align*}
   \frac{T^\epsilon}{\sqrt{NM}} \sum_{\substack{b,d,r \le T^{2} \\ al | dr \\ v|abl}}(abl)^{\frac{7}{32}}  \frac{lv^\half}{(dr)^\half} X^2 JK \ll  \sum_{\substack{b,d,r \le T^{2} \\ al | dr \\ v|abl}}(abl)^{\frac{7}{32}} \frac{(NM)^\half}{T^{2-\epsilon} a^2b^2l (dr)^\frac{3}{2} v^\half}.
   \end{align*}
This is the same as (\ref{c4}) and so we arrive at the same conclusions (\ref{summary1}-\ref{summary2}).

{\bf Summary.} We have reduced everything to proving that the following sum is bounded by a negative power of $T$:
 \begin{align*}
\frac{T}{\sqrt{NM}}  \sum_{\substack{K < \kk < 2K \\  }} \hat{\Phi}(t_g) A_f(\kk ,1)   \frac{\sqrt{\jj\kk drv}}{\cosh(\pi t_g)} \rho_g(\kk v)\overline{\rho_g(\jj dr)}
 \end{align*}
for any Maass cusp form $g$ of level $q\ll T^\epsilon$ and spectral parameter $t_g\ll T^\epsilon$, and any positive integers $T^{2-\epsilon} < K < T^{2+\epsilon}$, $\jj, v, d, r\ll T^{\epsilon}$. In the new ranges we have (remember that $\hat{\Phi}(t_g)$ depends implicitly on $\kk$) that
\begin{align*}
\frac{d}{d\kk} \hat{\Phi}(t_g)  \ll T^{-1+\epsilon}  \kk^{-1},
\end{align*}
by (\ref{jkders}). Thus by partial summation it suffices to show that 
 \begin{align}
\label{coffsum} T^{-2+\epsilon}  \sum_{K < \kk < 2K } A_f(\kk ,1)   \frac{\sqrt{\jj\kk drv}}{\cosh(\pi t_g)} \rho_g(\kk v)\overline{\rho_g(\jj dr)} \ll T^{-\delta}
 \end{align}
 for some $\delta>0$. 

\subsection{Generalized Lindel\"{o}f Hypothesis} \label{glhsection}

In order to use $L$-functions to obtain the required cancellation in (\ref{coffsum}), we must work with primitive Hecke cusp forms. To this end, we can take $g$ to be an element of the special basis described in Lemma \ref{bases}. That is, $g=h_c$ for some newform $h$ of level dividing $q$. Writing out $h_c$ has a linear combination of $h|_b$ as in the lemma, and using the fact the coefficients in this linear sum are small, it suffices to prove that
 \begin{align}
 T^{-2+\epsilon}  \sum_{K < \kk < 2K } A_f(\kk ,1)   \frac{\sqrt{\jj\kk drv}}{\cosh(\pi t_{h})} \rho_{h|_b}(\kk v)\overline{\rho_{h|_b}(\jj dr)} \ll T^{-\delta}
 \end{align}
for some newform $h$ of level $q'|q$ (so that $q'<T^{\epsilon}$) and spectral parameter $t_{h}=t_g<T^{\epsilon}$, and some positive integer $b<T^{\epsilon}$. If $\rho_h(n)$ are the Fourier coefficients of $h$, then the Fourier coefficients of $h|_b$ (see definition (\ref{slasdef})) are $\rho_h(\frac{n}{b})$ for $b|n$ and $0$ otherwise. Thus it suffices to prove that
 \begin{align*}
T^{-2+\epsilon}  \sum_{K < \kk < 2K }  A_f(b\kk ,1)   \frac{\sqrt{\jj\kk drv}}{\cosh(\pi t_h)} \rho_h(\kk v)\overline{\rho_h(\jj dr)} \ll T^{-\delta}
 \end{align*}
 for some $\delta>0$. We have 
 \begin{align*}
 \frac{1}{\cosh(\pi t_g)}  \sqrt{\kk v}\rho_h(\kk v) \sqrt{\jj dr} \ \overline{\rho_h(\jj dr)}=  \frac{1}{\cosh(\pi t_h)} |\rho_h(1)|^2  \lambda_h(\kk v)   \overline{\lambda_h(\jj dr)},
 \end{align*}
 where $\lambda_h(n)$ are the Hecke eigenvalues of $h$. By the standard (Rankin-Selberg) bound
 \begin{align*}
 \frac{|\rho_h(1)|^2 }{\cosh(\pi t_h)} \ll T^\epsilon,
 \end{align*}
we find that it suffices to prove that
 \begin{align*}
T^{-2+\epsilon}  \sum_{K < \kk < 2K }  A_f(b\kk ,1)   \lambda_h(\kk v) \ll T^{-\delta}
 \end{align*}
 for some $\delta>0$. Thus the trivial bound is $O(T^\epsilon)$, and to obtain further cancellation it suffices to partition the interval $[1,2]$ using smooth bump functions $\psi$ and bound the following by a negative power of $T$: 
\begin{align}
\nonumber &T^{-2+\epsilon}  \sum_{\kk \ge 1 } \psi\Big(\frac{\kk}{K}\Big)  A_f(b\kk ,1)   \lambda_h(\kk v) \\ 
\label{finalint} &= T^{-2+\epsilon}   \frac{1}{2\pi i} \int_{(2)}  \tilde{\psi}(s) K^s D(s) \frac{ds}{s},
 \end{align}
 where $\tilde{\psi}$ denotes the Mellin transform of $\psi$ and 
 \begin{align*}
 D(s)= \sum_{\kk \ge 1} \frac{A_f(b\kk ,1)   \lambda_h(\kk v)}{\kk^s}.
 \end{align*}
 For any integers $n,m$, let $n|m^\infty$ mean that $p|n\Rightarrow p|m$. By Hecke multiplicativity, we have
 \begin{align*}
 D(s) = \sum_{\substack{\kk \ge 1\\ (\kk,bv)=1}} \frac{A_f(\kk,1)\lambda_h(\kk)}{\kk^s} \sum_{\substack{\kk \ge 1\\ \kk|(bv)^\infty}} \frac{A_f(b\kk ,1)   \lambda_h(\kk v)}{\kk^s}.
 \end{align*}
 Comparing Euler products on both sides, we may write
 \begin{align*}
\sum_{\substack{\kk \ge 1\\ (\kk,bv)=1}} \frac{A_f(\kk,1)\lambda_h(\kk)}{\kk^s} = L(s,\sym^2 f\times h) G_1(s),
 \end{align*}
 where $G_1(s)$ is a Dirichlet series which absolutely converges for $\Re s\ge 1-\delta$ for some $\delta>0$ and satisfies $G_1(s)\ll 1$ in this half-plane. For any integer $n$ and prime $p$, let $n_p$ denote the $p$-part of $n$. That is, $n_p|n$ and $n_p p\nmid n$. We write
  \begin{align*}
G_2(s) = \prod_{p|bv} \Big( \frac{A(b_p,1)\lambda_h(v_p)}{1} + \frac{A(b_pp,1)\lambda_h(v_pp)}{p^s} +\frac{A(b_pp^2,1)\lambda_h(v_pp^2)}{p^{2s}} +\ldots\Big)
 \end{align*}
 and note that for $\Re(s)\ge \half$, we have 
 \begin{align*}
 G_2(s)\ll b^{\frac{7}{32}+\epsilon}v^{\frac{7}{64}+\epsilon} \exp(\sum_{p< \log bv} p^{-\half+\frac{7}{32}+\frac{7}{64}+\epsilon}) \ll T^\epsilon.
 \end{align*}
  Thus
 \begin{align*}
 D(s) = L(s,\sym^2 f\times h) G(s),
 \end{align*}
 where $G(s)=G_1(s)G_2(s)$ is a Dirichlet series which absolutely converges for $\Re s\ge 1-\delta$ for some $\delta>0$ and satisfies $G(s)\ll T^{\epsilon}$ in this half-plane. Thus moving the line of integration in (\ref{finalint}) to $\Re(s)=1-\delta$ and applying the GLH bound $L(s,\sym^2 f\times h) \ll (1+|s|)^\epsilon T^\epsilon$ there completes the proof (we save an absolute power of $T$, and all $\epsilon$'s can be adjusted). Note that what is actually required here is any subconvexity bound for $L(s,\sym^2 f\times h)$.
 
\subsection*{ Acknowledgment.} We thank Gergely Harcos for his encouragement to work on this problem.

\bibliographystyle{amsplain}
\bibliography{rwc}

\end{document}